\documentclass[10pt]{amsart}

\usepackage{amssymb,amsmath,amscd,amsfonts,a4,psfrag,graphicx,
latexsym,epsfig,color}
\usepackage[active]{srcltx}
\usepackage[all]{xy}

\everymath{\displaystyle}

\newtheorem{teo}{Theorem}[section]
\newtheorem{lem}[teo]{Lemma}
\newtheorem{cor}[teo]{Corollary}

\newtheorem{prop}[teo]{Proposition}
\newtheorem{defi}[teo]{Definition}

\newtheorem{remark}[teo]{Remark}
\newtheorem{remarks}[teo]{Remarks}

\newcommand{\fleche}[4]{                     
            \begin{array}{ccll} #1 & \longrightarrow & #2 \\   %
                         #3 &\longmapsto & #4          %
            \end{array}}

\newcommand{\cvd}{\hfill$\Box$}

\newcommand{\mc}{\mathbb{C}}
\newcommand{\mz}{\mathbb{Z}}
\newcommand{\mh}{\mathbb{H}}

\newcommand{\R}{\mathbb{R}}
\newcommand{\C}{\mathbb{C}}
\newcommand{\Z}{\mathbb{Z}}

\newcommand{\bw}{{\bf w}}

\newcommand{\bx}{{\bf x}}
\newcommand{\by}{{\bf y}}

\newcommand{\bW}{{\bf W}}

\newcommand{\Hh}{{\mathcal H}}

\newcommand{\Ll}{{\mathcal L}}

\newcommand{\Tt}{{\mathcal T}}

\newcommand{\rr}{{\bf r}}

\newcommand{\TG}{{\mathfrak T}}

\newcommand{\rG}{{\mathfrak r}}

\newcommand{\ra}{\rightarrow}
\newcommand{\Dim}{{\it Proof.\ }}

\def\cvd{\hfill$\Box$}

\title{On the quantum Teichm\"uller invariants of
  fibred cusped $3$-manifolds}

\author{St\'ephane Baseilhac$^1$, Riccardo Benedetti$^2$}

\begin{document}

\maketitle

\noindent $^1$ Institut Montpelli\'erain Alexander Grothendieck,
CNRS, Universit\'e de Montpellier 

(sbaseilh@univ-montp2.fr)

\smallskip

\noindent $^2$ Dipartimento di Matematica, Universit\`a di Pisa, Largo
Bruno Pontecorvo 5, 56127 Pisa, Italy (benedett@dm.unipi.it)


\begin{abstract}
  We show that the reduced quantum hyperbolic invariants of
  pseudo-Anosov diffeomorphisms of punctured surfaces are intertwiners
  of local representations of the quantum Teichm\"uller spaces. We
  characterize them as the only intertwiners that satisfy certain
  natural cut-and-paste operations of topological quantum field theories
  and such that their traces define invariants of mapping tori.
\end{abstract}

\tableofcontents

\section{Introduction}

Let $M$ be a {\it fibred cusped $3$-manifold}, that is, an oriented
non compact finite volume complete hyperbolic $3$-manifold which
fibres over the circle $S^1$.
In this paper we describe the relationships between
the quantum hyperbolic invariants of $M$ (QHI) (\cite{GT,AGT,NA}), and
the intertwiners of finite dimensional representations of the quantum
Teichm\"uller space of $S$ (\cite{B-L,B-B-L},
\cite{Filippo}), where $S$ is a fibre of some fibration of $M$ over
$S^1$.  Recall that $M$ is the interior of a compact manifold
$\bar M$ with boundary made by tori, and that $S$ is the interior
of a compact surface with boundary $\bar S$, properly embedded
in $\bar M$, so that $\bar S \cap \bar M$ is a family of `meridian'
curves on $\partial \bar M$.
\smallskip

The existence of such relationships has been expected for a long time,
and their precise formulation often thought as depending only on the
solution of a handful of technical problems. Let us recall the main
arguments underlying this opinion. Let $q$ be a primitive $N$-th root
of unity, where $N$ is odd. Then:
\begin{itemize}
\item The quantum Teichm\"uller space $\TG_{S}^q$ of $S$ consists
  of the set of Chekhov-Fock algebras $\TG_{\lambda}^q$, considered up 
  to a suitable equivalence relation (see Section \ref{QREP}), where $\lambda$ varies 
  among the ideal triangulations of $S$. A representation $\rho$ of
  $\TG_{S}^q$ is a suitable family of `compatible' representations
  $\rho = \{\rho_\lambda:\TG_{\lambda}^q\to {\rm End}(V_\lambda)\}_\lambda$.
 The irreducible
  representations of $\TG_{S}^q$ were classified up to
  isomorphism in \cite{B-L}, and so-called {\it local representations}
  were introduced and classified in \cite{B-B-L}. 
 \item The local representations are reducible, but more suited to
  perform the usual cut-and-paste operations of topological quantum field theories. In this
  framework, it is natural to restrict to a distinguished family of intertwiners,
satisfying some functorial
  properties with respect to the inclusion of triangulated
  subsurfaces. Theorem 20 of \cite{B-B-L} proposed such a family, and
  argued that for every couple $(\rho,\rho')$ of isomorphic local
  representations of  $\TG_{S}^q$ and every couple $(\lambda,\lambda')$
  of ideal triangulations of $S$, there is a unique projective class of
  intertwiners in the
  family which are defined from the representation space $V_\lambda$ of $\rho_\lambda$ to the representation space $V_{\lambda'}$ of $\rho'_{\lambda'}$. 
\smallskip

  Let us call the elements of these classes {\it QT intertwiners}.
  \smallskip
  
\item The diffeomorphisms of $S$ act on the local representations of
  $\TG_{S}^q$. Denote by $\phi$ the monodromy of some fibration of $M$ over $S^1$.
  If a local representation
  $\rho$ is isomorphic to its pull-back $\phi^*\rho$ under $\phi$ (an
  instance is canonically
  associated to the hyperbolic holonomy of $M$ \cite{B-L,B-B-L}), 
  then the trace of a suitably normalized QT intertwiner which relates 
   $\rho_{\lambda}$ to $\phi^*\rho_{\lambda'}$, $\lambda':=\phi(\lambda)$, 
 should be an invariant of $M$ (possibly
  depending on $\phi$) and the isomorphism class of $\rho$.
\item The QT intertwiners can be constructed by composing elementary intertwiners, 
associated to the diagonal exchanges in a sequence relating two ideal
  triangulations of $S$, and these were identified with the Kashaev
  $6j$-symbols in \cite{Ba}.
\item The QHI of $M$ can be defined by {\it QH state sums} over ``layered''
  triangulations of $M$, that is, roughly, $3$-dimensional ideal triangulations that
  realize sequences of diagonal exchanges relating the surface ideal triangulations
  $\lambda$ and $\phi(\lambda)$, in such a way that every diagonal exchange is associated to a tetrahedron of the triangulation. The main building
  blocks of the QH state sums are tensors called {\it matrix
    dilogarithms}, carried by the (suitably decorated) tetrahedra,
  which were derived in \cite{Top} from the Kashaev $6j$-symbols.
 \end{itemize}
 
In conclusion, the QHI of $M$ would coincide with the traces of QT
intertwiners of local representations of $\TG_{S}^q$ (at least when $M$ is equipped with
the hyperbolic holonomy). Only a suitable
choice of normalization of the QT intertwiners and an explicit correspondence between
formulas would be missing.

Another fact pointing in the same direction is that both 
the isomorphism classes of local representations of the quantum Teichm\"uller spaces and the QHI
depend on similar geometric structures. 
Every local representation $\rho$ of $\TG_{S}^q$ defines an {\it augmented 
character} of $\pi_1(S)$ in
$PSL(2,\mc)$ \cite{B-B-L}, called its {\it holonomy}, given by a system of {\it (exponential) shear-bend} coordinates on
every ideal triangulation $\lambda$ of $S$. Moreover $\rho$ has a {\it load} $h_\rho$,
defined as a $N$-th root of the product of the share-bend coordinates on any ideal triangulation $\lambda$ of $S$,
and which does not depend on the choice of $\lambda$.
The holonomy and the load classify $\rho$ up to isomorphism. 
 On another hand, 
the QHI of any cusped $3$-manifold $M$ depend on a choice of augmented character $\rG$ of $\pi_1(M)$ in
$PSL(2,\mc)$, lying in the
irreducible component of the variety of augmented characters containing the
hyperbolic holonomy. If $S$ is a fibre of a fibration of $M$
over $S^1$, $\rG$ can be recoved from its restriction to $\pi_1(S)$.
\smallskip

{\it However, other facts suggested that the relations between the QT intertwiners and the QHI would be 
subtler than expected}:

\begin{itemize}
\item 
The matrix dilogarithms are obtained from the Kashaev
$6j$-symbols by rewriting them in a non trivial way as tensors
depending on geometrically meaningful parameters (which are $N$-th
roots of hyperbolic shape parameters), {\it but also} by applying a
{\it symmetrization} procedure. This procedure is necessary to get the
invariance of the QH state sums over any cusped $3$-manifold. 
It involves several additional
choices, and eventually forces to multiply the QH state sums by a
normalization factor. It implies furthermore that the QHI must actually
depend not only on the choice of an augmented character $\rG$ of
$\pi_1(M)$ in $PSL(2,\mc)$, {\it but also} on some cohomological
classes, called {\it weights}, which cannot be recovered from the
holonomy of local representations of the quantum Teichm\"uller
spaces. The weights gauge out the additional choices. 
Surprisingly,
both the normalization factor of the QH state sums and the
weights did not appear in quantum Teichm\"uller theory.

\item An issue about the uniqueness of the projective classes
of QT intertwiners has been fixed recently. By a careful analysis
of the definition of local representations, and adapting the
construction of the QT intertwiners from \cite{B-B-L}, for every couple
$(\rho_\lambda, \rho'_{\lambda'})$ as above, it is constructed in
\cite{Filippo} a unique, minimal set
$\mathcal{L}_{\lambda\lambda'}^{\rho\rho'}$ of projective classes of
QT intertwiners which relate $\rho_\lambda$ to $\rho_{\lambda'}'$, and a free
transitive action
\begin{equation}\label{action1}
  \psi_{\lambda\lambda'}^{\rho\rho'}\colon H_1(S;\mz/N\mz) \times
  \mathcal{L}_{\lambda\lambda'}^{\rho\rho'} \to \mathcal{L}_{\lambda\lambda'}^{\rho\rho'}.
\end{equation}
In particular, the set $\mathcal{L}_{\lambda\lambda'}^{\rho\rho'}$ is far from being reduced to one point. 
\end{itemize}

\smallskip

The meaning of the normalization factor of the QH state sums has been 
better understood in \cite{NA}.  In that paper we showed that for every cusped $3$-manifold
$M$, the QH state sums of $M$ {\it without the
normalization factor} define finer invariants $\Hh_N^{red}(M,\rG,\kappa;\mathfrak{s})$,
called {\it reduced QH invariants}, well-defined up to multiplication
by $4N$-th roots of unity. They depend on an augmented
$PSL(2,\mc)$-character $\rG$ of $\pi_1(M)$, a class $\kappa\in
H^1(\partial \bar M;\mc^*)$, $\mc^*=\mc \setminus \{0\}$ being the multiplicative
group, such that $\kappa^N$ coincides up to sign
with the restriction of $\rG$ to $\pi_1(\partial \bar M)$, and an additional structure $\mathfrak{s}$
on $M$, called a {\it non ambiguous structure}. 

When $M$ fibres over $S^1$, a natural choice of $\mathfrak{s}$
is a fibration of $M$, say with monodromy $\phi$. Let us denote by
$M_\phi$ the realization of $M$ as the mapping torus of $\phi$,
equipped with the corresponding fibration, and the reduced QH invariant
by $\Hh_N^{red}(M_\phi,\rG,\kappa)$. The latter depends directly on the monodromy $\phi$, 
and hence is apparently more suited to investigate the relationships with the QT intertwiners. 
Note, however, that the class $\kappa$ is residual of the cohomological weights involved in the
definition of the unreduced QHI, and so the above considerations about the weights still apply.

Denote by $S$ a fibre of $M_\phi$.
Rather than considering $M_\phi$ one
can consider as well the mapping cylinder $C_\phi$ of $\phi$. In this
case, instead of a scalar invariant there is a {\it (reduced) QH-operator}
$$\Hh_N^{red}(T_{C_\phi},\tilde b,\bw)\in {\rm End}\left((\C^N)^{\otimes 2m}\right),
\quad m:=-\chi(S),$$
which is such that
\begin{equation}\label{traceform}
\Hh_N^{red}(M_{\phi},\rG,\kappa) = {\rm Trace}\left(\Hh_N^{red}(T_{C_\phi},\tilde b,\bw)\right).
\end{equation}
More precisely, in the spirit of topological quantum field theories, for every odd $N\geq 3$ the QHI define a contravariant functor from the category of (decorated) $(2 + 1)$-cobordisms to the category of vector spaces. In particular, it associates vector spaces $U_\lambda$ and $U_{\lambda'}$ to the triangulated surfaces $(S,\lambda)$ and $(S,\lambda')$ respectively, where $\lambda':=\phi(\lambda)$, and a QH-operator $\Hh_N^{red}(T_{C_\phi},\tilde b,\bw)\colon U_{\lambda'} \ra U_\lambda$ to the cylinder $C_\phi$. The diffeomorphism $\phi$ gives an identification between $U_\lambda$ and $U_{\lambda'}$, so one can consider $\Hh_N^{red}(T_{C_\phi},\tilde b,\bw)$ as an element of ${\rm End}\left(U_\lambda\right)$. Moreover, the construction of the functor provides a natural identification of $U_\lambda$ with $(\C^N)^{\otimes 2m}$. 

The operator $\Hh_N^{red}(T_{C_\phi},\tilde b,\bw)$
is defined by means of QH state sums over layered {\it QH-triangulations}
$(T_{C_\phi},\tilde b,\bw)$ of $C_\phi$. Here, by `layered QH-triangulation' we mean that: $T_{C_\phi}$ is a triangulation of $C_\phi$ that produces a layered triangulation of $M_\phi$ under the quotient map $C_\phi\rightarrow M_\phi$; $\bw$ is a suitable system of
$N$-th roots of shape parameters, solving the Thurston edge equations
associated to the triangulation of $M_\phi$ induced by $T_{C_\phi}$,
and which actually encodes both $\rG$ and $\kappa$;
finally $\tilde b$ is a way of ordering the elements of $\bw$ in each
tetrahedron of $T_{C_\phi}$, called {\it
  weak-branching}. The complete definition of $(T_{C_\phi},\tilde b,\bw)$ is given in Section \ref{quantum-hyp}.

\smallskip

Our results describe the relationships between the QH-operators
$\Hh_N^{red}(T_{C_\phi},\tilde b,\bw)$ and the set $\mathcal{L}_{\lambda\lambda'}^{\rho\rho'}$
of QT intertwiners. First note that these are assigned to $(\lambda,\lambda')$ in a {\it covariant} way, as they go from the representation space $V_\lambda$ of $\rho_\lambda$ to the representation space $V_{\lambda'}$ of $\rho'_{\lambda'}$. For every finite dimensional vector space
$V$, denote by $V'$ its dual space. The natural way to convert a QH operator $\Hh_N^{red}(T_{C_\phi},\tilde b,\bw)$
into a covariant one (preserving the value of its trace) is to consider its {\it dual}, or
{\it transposed operator}  
$$\Hh_N^{red}(T_{C_\phi},\tilde b,\bw)^T \in {\rm End}\left(U_\lambda'\right).$$
Using the definitions of $U_\lambda'$ in the QHI and of $V_\lambda$ in the theory of local representation of $\TG_{\lambda}^q$ one gets natural identifications $U_\lambda' \cong \left((\C^N)'\right)^{\otimes 2m} \cong V_\lambda$. Under all these identifications 
$$\Hh_N^{red}(T_{C_\phi},\tilde b,\bw)^T\in {\rm End}\left(\left((\C^N)'\right)^{\otimes 2m}\right)$$
is a QT intertwiner. Clearly, the matrix elements of these tensors
in the respective canonical bases of $\left((\C^N)'\right)^{\otimes 2m}$
and $(\C^N)^{\otimes 2m}$ are related by
\begin{equation}\label{transposed}
(\Hh_N^{red}(T_{C_\phi},\tilde b,\bw)^T)^{i_1\ldots i_{2m}}_{j_1\ldots j_{2m}} =
\Hh_N^{red}(T_{C_\phi},\tilde b,\bw)_{i_1\ldots i_{2m}}^{j_1\ldots j_{2m}}.
\end{equation}

We can now state our results precisely. Let $S$ be a punctured surface of negative Euler characteristic,
$\phi$ a pseudo-Anosov diffeomorphism of $S$, $M_\phi$ the mapping
torus of $\phi$ equipped with its fibration over $S^1$ with monodromy
$\phi$, and consider triples $(M_\phi,\rG,\kappa)$ and
$(T_{C_\phi},\tilde b,\bw)$ as above. Denote by $\lambda$, $\lambda'$
the ideal triangulations of $S$ given by the restriction to $S\times
\{0\}$ and $\phi(S)\times \{1\}$ of the (layered) triangulation
$T_{C_\phi}$. 
\begin{teo}\label{MAINTEO}
  (1) The QH-triangulation $(T_{C_\phi},\tilde b,\bw)$ determines
  representations $\rho_\lambda$ and $\rho_{\lambda'}$ of $\TG_\lambda^q$ and $\TG_{\lambda'}^q$ respectively, belonging to a
  local representation $\rho$  of $\TG_{S}^q$, such that $\rho_{\lambda}$ is isomorphic to $\phi^*\rho_{\lambda'}$, and acting on $V_\lambda= V_{\lambda'}=(\left((\C^N)'\right)^{\otimes 2m}$.  
  Moreover the transposed operator $\Hh_N^{red}(T_{C_\phi},\tilde b,\bw)^T$, considered as an element of ${\rm Hom}(V_\lambda,V_{\lambda'})$, is a QT intertwiner which intertwins the representations $\rho_\lambda$ and $\rho_{\lambda'}$.

(2) For any other choice of weak branching $\tilde b'$, the operator $\Hh_N^{red}(T_{C_\phi},\tilde b',\bw)^T$ intertwins local
  representations canonically isomorphic to $\rho_\lambda$,
  $\rho_{\lambda'}$ respectively.
\end{teo} 
We call $\Hh^\rho_{\lambda,\lambda'}:=\Hh_N^{red}(T_{C_\phi},\tilde b,\bw)^T$
a \emph{QHI intertwiner}.
\smallskip

Denote by $X_0(M)$ the (unique) irreducible component of the variety
of augmented $PSL(2,\mc)$-characters of $M$ that contains the discrete faithful
holonomy $\rG_h$. Denote by $X(S)$ the variety of augmented characters of $S$, and
by $i^*\colon X_0(M) \rightarrow X(S)$ the restriction map. Recall the residual weight $\kappa\in
H^1(\partial \bar M_\phi;\mc^*)$ involved in the definition of the reduced QHI.

\begin{teo} \label{MAINTEO2}
  (1) There is a neighborhood of $i^*(\rG_h)$ in $i^*(X_0(M)) \subset
  X(S)$ such that, for any isomorphism class of local representations of $\TG_S^q$
  whose holonomy lies in this neighborhood, there is a representative
  $\rho$ of the class, representations $\rho_\lambda$, $\rho_{\lambda'}$ belonging to $\rho$, and a QHI intertwiner
  $\Hh_N^{red}(T_{C_\phi},\tilde b,\bw)^T$ intertwining $\rho_\lambda$
  and $\rho_{\lambda'}$ as above. The load of $\rho$ is determined
  by the values of the weight $\kappa$ at the meridian
  curves that form $\bar S \cap \bar M$.

(2) The set of QHI intertwiners $\mathcal{H}^{\rho}_{\lambda,\lambda'}$
is the subset of QT intertwiners which intertwin the representation $\rho_\lambda$ to $\rho_{\lambda'}$, and whose
    traces are well defined invariants of the triples $(M_\phi,\rG,\kappa)$
    such that the restriction of $\rG$ to $\pi_1(S)$ is the holonomy of
    $\rho$.
\end{teo}

Let us make a few comments on Theorem \ref{MAINTEO} and \ref{MAINTEO2}.
\smallskip

$\bullet$ A key point of Theorem \ref{MAINTEO} (1) is to realize the transposed of the matrix dilogarithms as intertwiners between local representations of the two Chekhov-Fock algebras that one can associate to an ideal square. To this aim we consider an operator theoretic formulation of the matrix dilogarithms. As these are related to the Kashaev $6j$-symbols (in a non trivial way), one could alternatively develop a proof based on the result of Bai \cite{Ba} relating the Kashaev $6j$-symbols with intertwiners as above. However, it is not immediate. Bai's result deals with local representations up to isomorphism, and does not provide explicit relations between actual representatives. Also the relation between the Kashaev $6j$-symbols and the intertwiners of local representations is expressed in abstract terms, using the cyclic representations of the Weyl algebra, and not in terms of the eventual geometrically relevant parameters, the $q$-shape parameters (see Section \ref{quantum-hyp}). 
 
$\bullet$ Theorem \ref{MAINTEO2} (1) holds true more generally for any isomorphism class of local representations whose holonomy lies in $i^*(B)$, where $B$ is a determined Zariski open subset of the {\it eigenvalue} subvariety $E(M)$ of $X_0(M)$, introduced in \cite{KT}. Let us note here that $E(M)=X_0(M)$ if $M$ has a single cusp, and in general, $\rG_h \in E(M)$, and dim$_\mc E(M)$ equals the number of cusps of $M$. For simplicity, in this paper we restrict to characters in a neighborhood of $i^*(\rG_h)$ in $i^*(X_0(M))$. The general case is easily deduced from the results of \cite{AGT}.

$\bullet$ We can  say that {\it the invariance property of the QHI selects preferred elements in the set of all traces of QT intertwiners} (which has no a priori basepoints). Note that, the action \eqref{action1} being transitive, it does not stabilizes the set of QHI intertwiners. 
As well as the reduced QHI of $M_\phi$ give rise to the QT intertwiners $\Hh_N^{red}(T_{C_\phi},\tilde b,\bw)^T$,  with the mild ambiguity by $4N$-th roots of unity factors,
the unreduced QHI of $M$ give rise to the QT intertwiners $\Hh_N(T_{C_\phi},\tilde b,\bw)^T$, in the same projective classes as the ``reduced" ones.
The unreduced QHI of $M$, hence the associated QT intertwiners, do not depend on the choice of a fibration of $M$. However they depend on the full set
of cohomological weights, which dominates the classes $\kappa$. The ratio unreduced/reduced QHI is a simpler invariant called {\it symmetry defect}, which can be used to study the dependence of the reduced QHI with respect to the fibration of $M_\phi$ (\cite{NA}). Finally, the ambiguity of the invariants up to multiplication by $4N$-th roots of unity may not be sharp (see \cite{AGT} for some  improvements in the case of the unreduced QHI).
\medskip

The theorems above have the following consequences. 
\smallskip

Recall that the reduced QH invariants satisfy the identity $\Hh_N^{red}(M_{\phi},\rG,\kappa) = {\rm Trace}\big(\Hh_N^{red}(T_{C_\phi},\tilde b,\bw)\big)$, and that $M_\phi$ is the interior of a compact manifold $\bar M_\phi$ with boundary made by tori. Let us call {\it longitude} any simple closed curve in $\partial \bar M_\phi$ intersecting a fibre of $\bar M_\phi$ in exactly one point. We have:

\begin{cor}\label{cor0} The reduced QH invariants $\Hh_N^{red}(M_{\phi},\rG,\kappa)$ do not depend on the values of the weight $\kappa$ on longitudes.
\end{cor}

We can also use the QH-operators $\Hh_N^{red}(T_{C_\phi},\tilde b,\bw)$ to build finer invariants, associated to any irreducible representation of $\TG_S^q$. Namely, recall that every reducible representation is canonically the direct sum of its {\it isotypical} 
components (the maximal direct sums of isomorphic irreducible summands). 
Every intertwiner fixes globally each isotypical component. Then, we consider the \emph{isotypical intertwiners}
$$L^\phi_{\rho_\lambda(\mu)} := \Hh_N^{red}(T_{C_\phi},\tilde b,\bw)^T_{\vert \rho_\lambda(\mu)}$$
obtained by restriction to the isotypical components $\rho_\lambda(\mu)$ of $\rho_\lambda$, 
associated to irreducible representations $\mu$ of $\TG_{\lambda}^q$. Any 
irreducible representation of $\TG_\lambda^q$ can be embedded as a summand of some local 
representation (\cite{Tou}), and so has a corresponding isotypical component. 

We have:

\begin{cor}\label{cor}(1) The trace of $L^\phi_{\rho_\lambda(\mu)}$ is an invariant of $(M_{\phi},\rG,\kappa)$ and $\mu$, well-defined up to multiplication by $4N$-th roots of unity. It depends on the isotopy class of $\phi$ and satisfies
  $$\Hh_N^{red}(M_{\phi},\rG,\kappa)  = \sum_{\rho_\lambda(\mu)\subset \rho_\lambda}{\rm Trace}\left(L^\phi_{\rho_\lambda(\mu)}\right).$$

(2) The invariants ${\rm Trace}(L^\phi_{\rho_\lambda(\mu)})$ do not depend on the values of the weight $\kappa$ on longitudes. 
\end{cor}
\smallskip

In perspective, another application of Theorem \ref{MAINTEO} would be to study the representations of the Kauffman bracket skein algebras defined by means of the QHI intertwiners, by using the results of \cite{Tou}. 
\smallskip 

The background material is recalled in Section \ref{top-comb} and \ref{BACKRES}. The proofs are in Section \ref{PF}. 
\medskip

\noindent {\bf Remark about the parameter $q$.} In this paper we denote by $q$ an arbitrary primitive $N$-th root of unity, where $N\geq 3$ is odd. Our orientation conventions used to define the quantum Teichm\"uller space, in the relations \eqref{or1} and \eqref{or2}, imply that $q$ corresponds to $q^{-1}$ in some papers about quantum Teichm\"uller theory, eg. in \cite{B-L} or \cite{Filippo}.
Our choice of $q$ is motivated by Thurston's relations between shape parameters in hyperbolic geometry;
in the quantum case, this choice gives the most natural form to the tetrahedral and edge relations between $q$-shape parameters (see Section \ref{quantum-hyp}).

\section{Topological-combinatorial support}\label{top-comb}
We fix a compact closed oriented smooth surface $S_0$ of genus $g$,
and a subset $P=\{p_1, \dots, p_r\}\subset S_0$ of $r\geq 1$ marked points. We denote by $S$ the punctured surface S$_0\setminus P$. We assume that
$$m:=-\chi(S)= (2g - 2) + r>0.$$

A diffeomorphism $\phi_0:S_0\to S_0$ such that $\phi_0(p_j)=p_j$ for every $j$ induces a diffeomorphism $\phi: S\to S$.  Consider the cylinder $C_0:= S_0 \times [0,1]$ with the product orientation. Denote by $M_0$ the mapping torus of $\phi_0$ with the induced orientation. That is, $M_0:=C_0/\!\!\sim_{\phi_0}$, where $(x,0)\sim_{\phi_0} (y,1)$ if $y=\phi_0(x)$, $x, y\in S_0$. Let $\hat  M_0$ be the space obtained by collapsing to one point $x_j$ the image in
$  M_0$ of each line ${p_j}\times [0,1]$ in $C_0$. Then $\hat   M_0$ is a pseudo-manifold; $X=\{x_1,\dots ,x_r\}$ is the set of singular (ie. non manifold) points of $\hat   M_0$.
Every singular point $x_j$ has a conical neighbourhood homeomorphic to the quotient of $T^2_j \times [0,1]$ by the equivalence relation identifying $T^2_j \times \{0\}$ to a point, where $T^2_j$ is a $2$-torus. The point $x_j$ corresponds to the coset of $T^2_j \times \{0\}$, and $T^2_j\times \{1\}$ is called the {\it link} of $x_j$ in $\hat M_0$. So $M= \hat  M_0 \setminus X$ is the interior of a compact manifold $\bar M$ with boundary formed by $r$ tori. In fact $M$ is the mapping torus $M_\phi=C/\!\!\sim_\phi$, where $C=S\times [0,1]$. From now on, we denote by $M_\phi$ the manifold $M$ endowed with this fibration over the circle $S^1$, and by $C_\phi$ the mapping cylinder of $\phi$. 
\medskip

\begin{figure}[ht]
\begin{center}
\includegraphics[width=7cm]{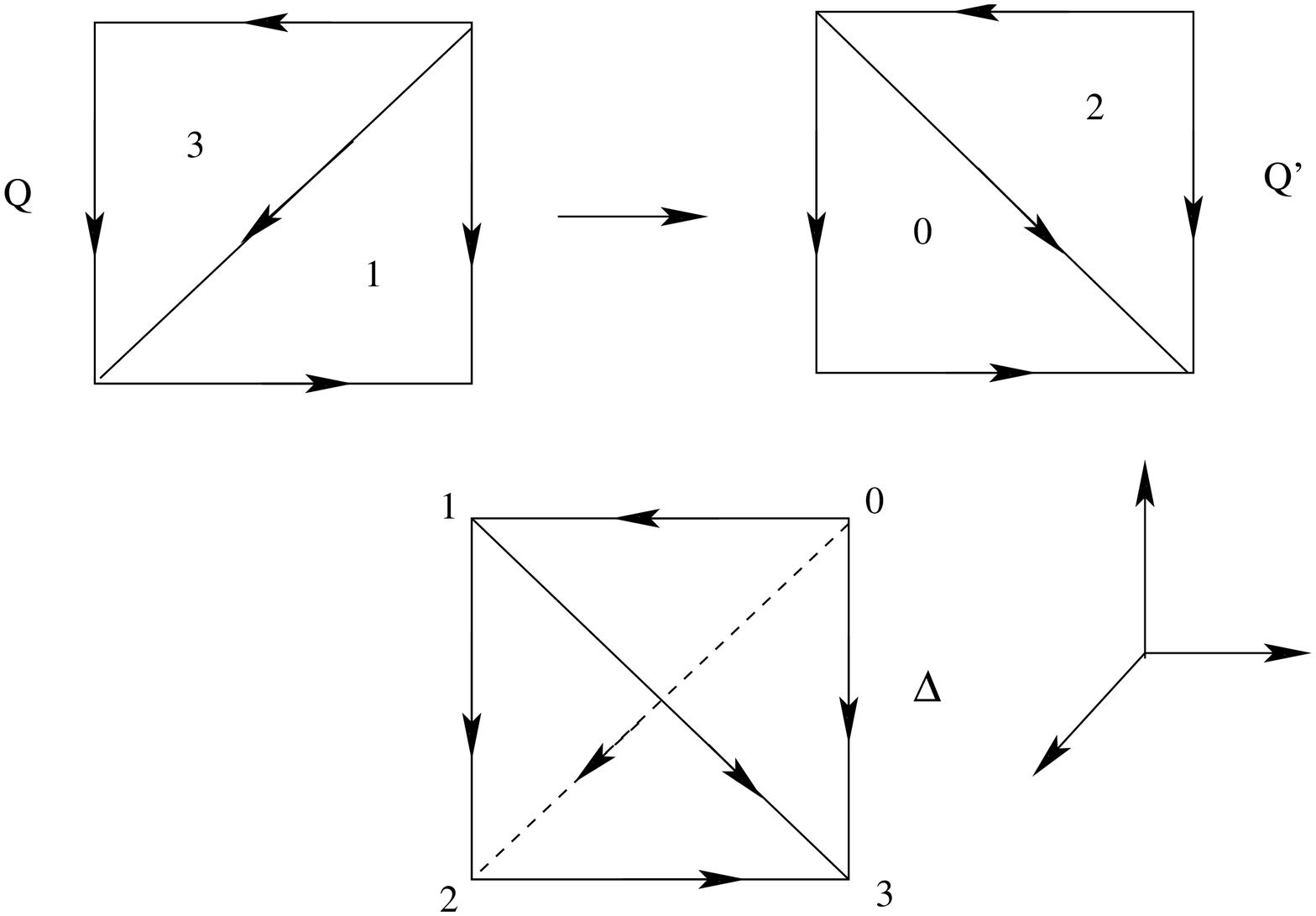}
\caption{\label{typefig} }
\end{center}
\end{figure}

Layered triangulations of $M_\phi$ can be constructed as follows (see for instance \cite{LA}, where one can find also a proof of the
well known fact that  two ideal triangulations of $S$ are connected by a finite chain of diagonal exchanges).
Consider the set of triangulations $\lambda$ of $S_0$ whose sets of vertices coincide
with $P$. By removing $P$, such a triangulation $\lambda$ is also said an 
{\it ideal triangulation} of $S$; it has $3m$ edges and $2m$ triangles. Denote by $\lambda'= \phi_0(\lambda)$ the image ideal triangulation of $S$.
Let us fix a chain of diagonal exchanges connecting $\lambda$ to $\lambda'$:
$$\lambda = \lambda_0 \to  \lambda_1 \to \dots \to \lambda_k= \lambda'.$$ 
Possibly by performing some additional diagonal exchanges followed by their inverses, we can assume that this chain is ``full'',
in the sense that every edge of $\lambda$ supports some diagonal exchange in the chain. Then the chain induces a $3$-dimensional
triangulation $T$ of $\hat M_0$, whose tetrahedra are obtained by ``superposing'', for every $j$, the two squares of $\lambda_j$, $\lambda_{j+1}$ involved in
the diagonal exchange $\lambda_j\rightarrow \lambda_{j+1}$ (see Figure \ref{typefig}). The set of vertices of $T$ coincides with $X$.
The ideal triangulation $T\setminus X$ is called a {\it layered triangulation} of $M_\phi$. It lifts to a layered triangulation $T_{C_\phi}$ of the cylinder $C_\phi$.
\medskip

For each tetrahedron of $T$, we call \emph{abstract tetrahedron} its
underlying simplicial set considered independently of the face
pairings in $T$. Similarly we call abstract edge or face of $T$ any
edge or face of an abstract tetrahedra of $T$. Given any edge $e$ of
$T$, we write $E \to e$ to mean that an abstract edge $E$ is
identified to $e$ in the triangulation $T$.
\medskip

From now on we assume that the diffeomorphism $\phi$ of $S$ is {\it pseudo-Anosov},
so that $M_\phi$ is a hyperbolic manifold with $r$ cusps.

\medskip

By construction every (ideally triangulated) surface $(S,\lambda_j)$
is embedded in $T$ and $T_{C_\phi}$. The union of the surfaces $(S,\lambda_j)$
forms the $2$-skeletons of $T$ and $T_{C_\phi}$. Every surface $(S,\lambda_j)$ divides $M_\phi$ locally,
and the given orientations of $S$ and $M_\phi$ determines the (local) {\it positive
  side} of $S$ in $M$. In terms of $T_{C_\phi}$, the positive side of
$(S,\lambda)=(S,\lambda_0)$ lifts to a collar of $S\times \{0\}$,
while the {\it negative} side of $(S,\lambda')=(S,\lambda_k)$ lifts to
a collar of $S\times \{1\}$.
\medskip

The triangulation $T$ is naturally endowed with a {\it taut pre-branching} $\omega$. Let us recall briefly this notion (see \cite{NA}, \cite{LA}). A taut pre-branching is a choice of transverse co-orientation for each $2$-face of $T$, so that for every tetrahedron, exactly two of its $2$-face co-orientations point inward/outward; moreover, it is required that there are exactly two abstract diagonal edges $E\to e$ for every edge $e$ of $T$. Here, we call \emph{diagonal edges} the two edges of an abstract tetrahedron whose adjacent $2$-face co-orientations point both inward or both outward. The taut pre-branching $\omega$ of $T$ is defined by the system of transverse $2$-face co-orientations dual to the orientation of the embedded surfaces $(S,\lambda_j)$. 

We fix also a {\it weak branching} $\tilde b$ {\it compatible with} $\omega$ (see \cite{AGT}). Recall that $\tilde{b}$ is a choice of vertex ordering $b$ for each abstract tetrahedron of $T$, called {\it (local) branching}, satisfying the following constraint. The abstract $2$-faces of the tetrahedron have an orientation induced by $b$, given by their vertex orderings up to even permutations. Then, it is required that these orientations match, under the $2$-face pairings, with the orientations dual to $\omega$ (note that we do not require that the abstract vertex orderings match: this defines the stronger notion of {\it global branching}). Similarly, $\tilde{b}$ gives an orientation to every tetrahedron; without loss of generality, we can assume that it agrees with the orientation of $M$. A weak branching $\tilde b_j$ is also defined on every surface triangulation $\lambda_j$, by giving every triangle $\tau$ in $\lambda_j$ the branching induced by the tetrahedron lying on the positive side of $\tau$.
\medskip

In Figure \ref{typefig} we show a typical {\it local} configuration
occurring in $T$ or $T_{C_\phi}$, that is, a tetrahedron $\Delta$, represented under an orientation preserving embedding in $\R^3$, built by gluing squares $Q$, $Q'$ along the boundary, carrying triangulations $\lambda$, $\lambda'$ related by a diagonal exchange. The orientation of $\Delta$ agrees with the standard orientation of
$\R^3$, and the $2$-face co-orientations that define the pre-branching $\omega$ are dual to the counter-clockwise orientation of the
four faces of $\Delta$. We have chosen one branching $b$ of $\Delta$ (any edge being oriented from the lowest to the biggest endpoint) among the two
which are compatible with $\omega$ and induce the given orientation of $\Delta$.
The $2$-faces are ordered as the opposite vertices are, accordingly with $b$. 
\smallskip

We advertize that Figure \ref{typefig} will be used to support all computations. In such a situation we will denote
by $e_0$ the edge $[v_0,v_1]$, by $e_1$ the edge $[v_1,v_2]$, and by $e_2$ the edge $[v_0,v_2]$. 

\section{Background results}\label{BACKRES}

\subsection{Classical hyperbolic geometry}\label{class-hyp} Let $M$ be any cusped hyperbolic $3$-manifold with $r$ cusps, $r\geq 1$. 
\smallskip

We denote by $X(M)$ the variety of {\it augmented $PSL(2,\mc)$-characters} of $M$. Let us recall the definition. Consider the set of pairs $(\rr,\{\xi_\Gamma\}_{\Gamma \in \Pi})$, where 
$\rr\colon \pi_1(M)\rightarrow PSL(2,\mc)$ is a group homomorphism, $\Pi$ is the set of peripheral subgroups of $\pi_1(M)$ (associated to the boundary tori $T^2_i$ of $\bar M$, for all choices of base points of $\pi_1(T^2_i)$ and $\pi_1(M)$ and paths between them), and for every $\Gamma\in \Pi$, $\xi_\Gamma\in \mc P^1$ is a fixed point of $\rr(\Gamma)$ such that the assignment $\Gamma \rightarrow \xi_\Gamma$ is equivariant with respect to the action of $\pi_1(M)$ on $\Pi$ by conjugation and on $\mc P^1$ via $\rr$ by Moebius transformations. The set $\{(\rr,\{\xi_\Gamma\}_{\Gamma \in \Pi})\}$ is a complex affine algebraic set $R(M)$, with an action of $PSL(2,\mc)$ defined on a pair $(\rr,\{\xi_\Gamma\}_\Gamma)$ by conjugation on $\rr$ and Moebius transformation on $\{\xi_\Gamma\}_\Gamma$. Then, $X(M)$ is the algebro-geometric quotient of $R(M)$ by $PSL(2,\mc)$, that is, the set of closed points of the ring of invariant functions $\mc[R(M)]^{PSL(2,\mc)}$.  
\smallskip

Similarly we denote by $X(S)$ the variety of augmented characters of a surface $S$. 
\smallskip

In the setting of Section \ref{top-comb}, if $M=M_\phi$ and $S$ is a fibre of the fibration, then
the inclusion map $i\colon S=S\times \{0\} \hookrightarrow M_\phi$ induces a regular (restriction) map
$$i^*\colon X(M_\phi) \rightarrow X(S).$$ 
Recall that $\pi_1(M_\phi)$ is a HNN-extension of $\pi_1(S)$. That is, given a generating set $\gamma_1,\ldots,\gamma_u$ of $\pi_1(S)$ satisfying the relations $r_1,\ldots,r_v$, there is an isomorphism between $\pi_1(M_\phi)$ and the group generated by $\gamma_1,\ldots,\gamma_u$ and an element $t$ satisfying the relations $r_1,\ldots,r_v$ and $t\alpha t^{-1} = \phi_*(\alpha)$ for all $\alpha\in \pi_1(S)$, where $\phi_* \colon \pi_1(S)\ra \pi_1(S)$ is the isomorphism induced by $\phi$. Therefore, every representation $\rr\colon \pi_1(S)\rightarrow PSL(2,\mc)$ that may be extended to $M_\phi$ is such that $\rr(\pi_1(S))$ and $\rr\circ \phi_*(\pi_1(S))$ are conjugate subgroups of $PSL(2,\mc)$. Any augmented character of $S$ that may be extended to $M_\phi$ is thus a fixed point of the
map $\Phi_*\colon X(S) \rightarrow X(S)$ induced by the map $\rr\mapsto \rr\circ \phi_*$ on representations, and $\overline{i^*(X(M_\phi))}$ is a subvariety of ${\rm Fix}(\Phi_*)$.
\smallskip

Again in the setting of Section \ref{top-comb}, let $(T,\tilde b)$ be a weakly branched layered triangulation of the fibred cusped manifold $M_\phi$. Let us recall a few facts about the {\it gluing variety} $G(T,\tilde{b})$, that is, the set of solutions of the Thurston gluing equations supported by $(T,\tilde{b})$. It is a complex affine algebraic set of dimension greater than or equal to $r$. The points of $G(T,\tilde{b})$ are classically called {\it systems of shape parameters}. Their coordinates, the shape parameters, are scalars in $\mc\setminus \{0,1\}$ associated to the abstract edges of $T$. The shape parameters of opposite edges of a tetrahedron are equal, and the cyclically ordered triple of shape parameters of a tetrahedron encodes an isometry class of hyperbolic ideal tetrahedra.  

The set $G(T,\tilde{b})$ is defined by the following two sets of equations. For every branched tetrahedron, set $w_j := w(E_j)$, $j=0,1, 2$  (see before Figure \ref{typefig} for the ordering of the edge $E_j$). For every edge $e$ of $T$, define the {\it total shape parameter} $W(e)$ as the product of the shape parameters $w(E)$, where $E \to e$. Then we have:
\begin{itemize}\label{crel}
\item {\it (Tetrahedral equation)} For every tetrahedron and $j\in\{0,1,2\}$, $w_{j+1}=(1-w_j)^{-1}$ cyclically; hence $w_0w_1w_2=-1$. 
\item {\it (Edge equation)} For every edge $e$ of $T$, $W(e)=1$.
\end{itemize}

A point $w\in G(T,\tilde{b})$ determines a pseudo-developing map $F_w:\tilde{M}_\phi \rightarrow \mh^3$, where $\tilde{M}_\phi$ is the universal cover of $M_\phi$, and $F_w$ is well-defined up to post-composition with an orientation-preserving isometry of $\mh^3$. The map $F_w$ sends homeomorphically the edges of $\tilde T$ to complete geodesics, and it satisfies $F_w (g \tilde x) = \rr_w(g) F_w (\tilde x)$ for all $\tilde x\in \tilde M_\phi$, $g\in \pi_1(M_\phi)$, where $\rr_w\colon \pi_1(M_\phi)\to PSL(2,\mc)$ is a homomorphism. So $w$ encodes the conjugacy class of $\rr_w$. In fact, $F_w$ determines also some $\rr_w$-equivariant set $\{\xi_\Gamma\}_{\Gamma \in \Pi}$ as above, so that eventually the map $w\mapsto \rr_w$ can be lifted to a regular ``holonomy'' map
$$hol\colon G(T,\tilde{b})\rightarrow X(M_\phi).$$
We have (see Proposition 4.6 of \cite{AGT} when $M_\phi$ has a single cusp, and Remark 1.4 of \cite{NA} for the general case):
\begin{prop}\label{gluing-var} There is a subvariety $A$ of $G(T,\tilde{b})$ of dimension equal to the number of cusps of $M$, and a point $w_h\in A$, such that $hol(w_h)$ is the hyperbolic holonomy of $M_\phi$, and $hol_{\vert A}$ is a homeomorphism from a Zariski open subset of $A$ containing $w_h$ onto its image. 
\end{prop}

In particular, every point $w\in A$ encodes an augmented character of $M_\phi$; the algebraic closure of $hol(A)$ is the {\it eigenvalue} subvariety $E(M_\phi)$ of $X_0(M_\phi)$, the irreducible component of $X(M_\phi)$ containing the discrete faithful holonomy $\rr_h$ (see \cite{KT}). If $M_\phi$ has a single cusp, then $E(M_\phi)=X_0(M_\phi)$; in general $E(M_\phi)$ contains $\rr_h$ and has complex dimension equal to the number of cusps of $M_\phi$. 
\medskip 

It follows from Proposition \ref{gluing-var} that a point of $i^*(X_0(M_\phi))\subset X(S)$ close enough to $i^*(\rr_h)$ is encoded by a point $w\in A$.
\medskip

From now on, we consider only systems of shape parameters $w$ lying in the subvariety $A$ of $G(T,\tilde{b})$. Denote by $(T,\tilde{b},w)$ the layered triangulation $T$ 
of $M_\phi$ endowed with the weak branching $\tilde{b}$ and the labelling of the abstract edges of $T$ by a system of shape parameters $w$. For every $j=0,\dots, k$, consider the ideally triangulated surfaces $(S,\lambda_j)$ embedded into $(T,\tilde b, w)$ as in Section \ref{top-comb}.  For every edge $e$ of
$\lambda_j$, define the {\it lateral shape parameter} $W^+_j(e)$ as the product of the
shape parameters of the abstract edges $E \to e$ carried by the tetrahedra lying on the {\it positive} side of $(S,\lambda_j)$.

\begin{lem}\label{shear-bend} For any edge $e$ of $\lambda_j$, the (exponential) shear-bend coordinate of $e$ defined by any pleated surface $F_{w\vert \tilde S}$ (ie. any lift $(\tilde{S},\tilde \lambda_j)$ of $(S,\lambda_j)$ to $\tilde M_\phi$, and any pseudo-developing map $F_{w\vert \tilde S}\colon \tilde S \rightarrow \mh^3$), coincides with the parameter $\textstyle x^j(e):= -W^+_j(e)$.
\end{lem}

The proof follows from the definitions (see eg. \cite{B-L}). Note that $x^0(e)=x^k(\phi_0(e))$ for every edge $e$ of $\lambda=\lambda_0$. Also, the shear-bend coordinate of any edge $\tilde e$ of $\tilde \lambda_j$ is eventually `attached' to the corresponding edge $e$ of $\lambda_j$, because for different choices of $\tilde S$ or $F_w$ the images of $F_{w\vert \tilde S}$ differ only by an hyperbolic isometry. 

\begin{remarks}{\rm The parameter $x^j(e)$ is the {\it opposite} of $W^+_j(e)$ because the (oriented) bending angle along an edge of $\lambda_j$ is traditionally measured by the external dihedral angle $\pi-\theta$ (see eg. \cite{B-L,B-B-L}), whereas the shape parameters use the internal dihedral angle $\theta$. So $\theta=0$ when two adjacent ideal triangles $F (\tau)$ and $F(\tau')$ coincide, for triangles $\tau$, $\tau'$ of $(\tilde{S},\tilde{\lambda}_j)$.}
\end{remarks}

\subsection{Quantum hyperbolic geometry}\label{quantum-hyp}
We keep the  setting of the previous section.\smallskip

Recall that $N\geq 3$ is an odd integer, and $q$ a primitive $N$-th root of unity.
\smallskip

We begin with a few qualitative, structural features of the 
{\it reduced quantum hyperbolic state sum} $\Hh^{red}_N(T,\tilde b,\bw)$ defined in \cite{AGT,NA}. It is a regular rational function
defined on a covering of the gluing variety $G(T,\tilde b)$. The
points of this covering over a point $w\in  G(T,\tilde b)$ are certain systems of $N$-th roots $\bw(E)$ of
the shape parameters $w(E)$, which label the abstract edges $E$
of $T$  and are called {\it
  quantum shape parameters}. Alike the ``classical" ones, opposite edges
  of an abstract tetrahedron are given the same quantum shape parameter.
  Moreover,  the quantum shape parameters verify the following 
  relations, which are ``quantum" counterparts of the defining equations of the
  gluing variety. For every branched tetrahedron of 
  $(T,\tilde b)$ put
$\bw_j := \bw(E_j)$ (with the usual edge ordering fixed before Figure
\ref{typefig}). For every edge $e$ of $(T,\tilde b,\bw)$, define the
    {\it total quantum shape parameter} $\bW(e)$ as the product of the
    quantum shape parameters of the abstract edges $E \to e$.  Then we
    have:
\begin{itemize}\label{qrel}
\item {\it (Tetrahedral relation)} For every tetrahedron, $\bw_0\bw_1\bw_2= -q$.
\item {\it (Edge relation)} For every edge $e$ of $T$, $\bW(e)=q^{2}$.  
\end{itemize}
For simplicity we will work with systems of quantum shape parameters
$\bw$ such that the corresponding systems of ``classical'' shape parameters $w$ belong to a simply connected open neighborhood $A_0\subset A\subset
G(T,\tilde b)$ of $w_h$ (the general case of systems $\bw$ with
arbitrary $w\in A$ is described in \cite{AGT}).  
Note that $A_0$ is chosen so that $\bw(E)$ varies continuously with
$w\in A_0$. 
\begin{remarks}\label{qpar} {\rm  (1) For every $z\in \C\setminus \{0\}$, denote by $\log(z)$ the determination of the logarithm which
has the imaginary part in $]-\pi,\pi]$. There is a $\mz$-labelling $d$ of the abstract edges
of $T$ such that, for every system $\bw$ of quantum shape parameters over a point $w\in A_0$, and every abstract edge
$E$, we have
\begin{equation}\label{qshapeconst}
\bw(E) = \exp\left(\frac{1}{N}\left(\log(w(E))+\pi i (N+1)d(E)\right)\right).
\end{equation}
The above relations verified by the quantum shape parameters
can be equivalently rephrased in terms of the $\mz$-labelling $d$.
For any point $w\in A_0$, a $\mz$-labelling $d$ satisfying these
relations defines a system of quantum shape parameters over $w$ by the formula
\eqref{qshapeconst}. 

(2) In \cite{GT,AGT,NA} we solved (mainly in terms of the $\mz$-labellings $d$) the existence problem of triples $(T,\tilde b,\bw)$ for any cusped
hyperbolic $3$-manifold in the special case $q=-\exp(-i\pi/N)$. 
The same method works for an arbitrary $q$, up to minor changes.}
\end{remarks}
 
The following result summarizes in a qualitative way the invariance properties of $\Hh^{red}_N(T,\tilde b,\bw)$.
Denote by $A_{0,N}$ the set of systems of quantum shape parameters
over $A_0$. Recall from Section \ref{top-comb} that $M_\phi$ can be considered as the interior of a compact manifold
$\bar M_\phi$ bounded by tori. Fix a basis $(l_1,m_1),
\ldots , (l_r,m_r)$ of $\pi_1(\partial \bar M_\phi)$, and use it to
identify $H^1(\partial \bar{M}_\phi;\mc^*)$ with $(\mc^*)^{2r}$. Any augmented character $[(\rr,\{\xi_\Gamma\}_\Gamma)]\in X(M_\phi)$ determines the square of one of the two (reciprocally inverse) eigenvalues of $\rr(\gamma)$, for any non trivial simple closed curve $\gamma$ on $\partial \bar M_\phi$ and representative $\rr$ of the character. Namely, by taking $\rr$ in its conjugacy class so that $\rr(\pi_1(\partial \bar M_\phi))$ fixes the point $\infty$ on $\mathbb{C}P^1$, $\rr([\gamma])$ acts on $\mathbb{C}$ as $w\mapsto \gamma_\rr w + b$, where $\gamma_\rr\in \mathbb{C}^*$ and $b \in \mathbb{C}$. The coefficient $\gamma_\rr$ is that squared eigenvalue selected by $[(\rr,\{\xi_\Gamma\}_\Gamma)]$.

We have (see \cite{AGT} and \cite{NA}):
\begin{teo} \label{SSum}
  (1) There exists a determined regular rational map $\kappa_N\colon
  A_{0,N}\ra (\mc^*)^{2r}$ such that the image $\kappa_N(A_{0,N})$ is the open
  subset of $(\mc^*)^{2r}$ made of all the classes $\kappa\in
  H^1(\partial \bar{M}_\phi;\mc^*)$ such that $\kappa(l_j)^N$ and
  $\kappa(m_j)^N$ are equal up to a sign respectively to the squared
  eigenvalues selected by $[(\rr,\{\xi_\Gamma\}_\Gamma)]$ at the curves $l_j$ and $m_j$,
  where $\rG:=[(\rr,\{\xi_\Gamma\}_\Gamma)] =hol(w)$ for some $w\in A_0$.

(2) Given $\rG\in hol(A_0)$ and $\kappa\in H^1(\partial
  \bar{M}_\phi;\mc^*)$ as above, the value of $\Hh^{red}_N(T,\tilde
  b,\bw)$ is independent, up to multiplication by $4N$-th roots of
  unity, of the choice of $(T,\tilde b,\bw)$ among all weakly branched
  layered triangulations of $M_\phi$ endowed with a system of quantum
  shape parameters $\bw$ such that $hol(w) =
  \rG$ and $\kappa_N(\bw) =\kappa$.
\end{teo}
We denote the resulting invariant by $\Hh_N^{red}(M_\phi,\rG,\kappa):=
\Hh^{red}_N(T,\tilde b,\bw)$. By the results of \cite{NA}, {\it it actually
depends on the fibration of $M_\phi$}.
\medskip

The map $\kappa_N$ in Theorem \ref{SSum} is a lift to $A_{0,N}$ of $j^*\circ hol_{\vert A_0}$, where $j^*$ is the
restriction map $X(M_\phi) \ra X(\partial M_\phi)$. 
Here is a concrete way to compute $\kappa_N(\bw)(\alpha)$, for $\alpha \in H_1(\partial \bar{M}_\phi;\mz)$. 
It is well known that the {\it truncated tetrahedra} of an ideal triangulation $T$ of
$M_\phi$ provide a cell decomposition of $\bar M_\phi$ which restricts to 
a triangulation $\partial T$ of $\partial \bar M_\phi$.
Every abstract vertex of $T$ corresponds to an abstract triangle of $\partial T$ (at which a tetrahedron
has been ``truncated''); every vertex of such a triangle is contained in one abstract edge of $T$.
Represent the class $\alpha$ by an oriented simple closed curve $a$ on $\partial \bar M_\phi$,
transverse to $\partial T$.
The intersection of $a$ with a triangle $F$ of $\partial T$ is a collection of oriented arcs.
We can assume that none of these arcs enters and exits $F$ by a same edge.
Then, each one turns around a vertex of $F$.
If $\Delta$ is a tetrahedron of $T$ and $F$ corresponds to a vertex of $\Delta$,
for every vertex $v$ of $F$ we denote by $E_v $
the edge of $\Delta$ containing $v$, and write $a \rightarrow E_v$ to mean that some subarcs of
$a$ turn around $v$.
We count them algebraically, by using the orientation of $a$: if there are $s_+$ (resp. $s_-$)
such subarcs whose
orientation is compatible with (respectively, opposite to) the
orientation of $\partial \bar M_\phi$ as viewed from $v$, then we set $ind(a,v):=s_+-s_-$. Then
\begin{equation}\label{cuspweightform2}
\kappa_N(\bw)(\alpha) = \prod_{a \ra E_v}  \bw(E_v)^{ind(a,v)}.
\end{equation}
In the case where $\alpha$ is the class of a positively oriented meridian curve $m_i$ of
$\partial \bar M_\phi$, this gives
\begin{equation}\label{cuspweightform}
\kappa_N(\bw)(\alpha) = \prod_{m_i \ra E_v}  \bw(E_v).
\end{equation}

Let us give now more details about the definition of the reduced QH state sum $\Hh^{red}_N(T,\tilde b,\bw)$. 
We are going to do it in terms of the QH-operator
associated to the cylinder $C_\phi$, already mentioned in the Introduction.

By cutting $(T,\tilde b,\bw)$ along one of the surfaces $(S,\lambda_j)$
(that is, a triangulated fibre of $M_\phi$),
we get a {\it QH-triangulation} $(T_{C_\phi},\tilde b,\bw)$ of the cylinder $C_\phi$
having as ``source" boundary component $(S,\lambda)$, and as ``target" boundary component $(S,\lambda')$, where $\lambda=\lambda_j$ and $\lambda'=\phi(\lambda_j)$. It carries in a {\it contravariant way} the QH-operator
$\Hh_N^{red}(T_{C_\phi},\tilde b,\bw)\in {\rm End}\left((\C^N)^{\otimes 2m}\right)$, which is such that 
$${\rm Trace}\left(\Hh_N^{red}(T_{C_\phi},\tilde b,\bw)\right)=\Hh^{red}_N(T,\tilde b,\bw).$$
We define the QH-operator $\Hh_N^{red}(T_{C_\phi},\tilde b,\bw)$ by means of elementary tetrahedral and $2$-face operators, as follows. 

We can regard the QH triangulation $(T,\tilde b,\bw)$ of $M_\phi$ as a  {\it network of abstract QH-tetrahedra} $(\Delta,b,\bw)$, 
with gluing data along  the $2$-faces.
To each QH-tetrahedron $(\Delta,b,\bw)$ we associate a linear isomorphism called {\it basic matrix dilogarithm}
(``basic'' refers to the fact that we have dropped a symmetrization factor from the matrix dilogarithms that enter the definition of the unreduced QHI),
$$\Ll_N(\Delta,b,\bw)\colon V_2\otimes V_0 \rightarrow V_3 \otimes V_1$$
where $V_j$ is a copy of $\mc^N$ associated to the $2$-face of $(\Delta,b)$ opposite to the vertex $v_j$. The basic matrix dilogarithms will be defined explicitly below. For the moment, note that $V_1$, $V_3$ correspond to the $2$-faces with pre-branching co-orientation pointing inside $\Delta$. Concerning the gluing data, every $2$-face $F$ of $T$ is obtained by gluing a pair of abstract $2$-faces. Denote by $F_s$ and $F_t$ the ``source'' and ``target'' $2$-face of the pair, with respect to the transverse co-orientation defined by the pre-branching $\omega_{\tilde b}$. Denote by $V_{F_s}$ and $V_{F_t}$ the copies of $\mc^N$ associated to $F_s$ and $F_t$. The identification $F_s\to F_t$ is given by an even permutation on three elements, which encodes the image of the vertices of $F_s$ in $F_t$. So it can be encoded by an element $r(F)\in \Z/3\Z$.
The triangulation  $T_{C_\phi}$ of the cylinder $C_\phi$ has $2m$ free $2$-faces at both the source and target boundary components, $(S,\lambda)$ and $(S,\lambda')$. Note that at every $2$-face of  
$(S,\lambda)$ the pre-branching orientation
points inside $C_\phi$, while it points outside at the $2$-faces of $(S,\lambda')$. For every gluing $F_s\to F_t$ occuring at an internal $2$-face $F$ of $T_{C_\phi}$, as well as to every $2$-face $F$ of $(S,\lambda')$, we associate an endomorphism (again in contravariant way)
\begin{equation}\label{Qop}
Q_N^{r(F)} : V_{F_t} \to V_{F_s}.
\end{equation}
By definition, the QH-operator $\Hh_N^{red}(T_{C_\phi},\tilde b, \bw)$ is the total contraction of the network of tensors $\{\Ll_N(\Delta,b,\bw)\}_{\Delta}$ and $\{Q_N^{r(F)}\}_F$. In formulas:
\begin{equation}\label{Ssumformula}
\Hh_N^{red}(T_{C_\phi},\tilde b, \bw)  = \sum_s \prod_\Delta \Ll_N(\Delta,b,\bw)_s \prod_F Q_{N,s}^{r(F)}
\end{equation}
where the sum ranges over all maps 
$s\colon\{$abstract $2$-faces of $T\} \cup \{$$2$-faces of $(S,\lambda')\} \to \{0,\ldots,N-1\}$
(the {\it states} of $(T_{C_\phi},\tilde b, \bw)$), and $\Ll_N(\Delta,b,\bw)_s$ and $Q_{N,s}$ denote the entries of the tensors $\Ll_N(\Delta,b,\bw)$ and $Q_{N}$ selected by $s$, when the tensors are written in the canonical basis $\{e_j\}$ of $\mc^N$.  Note that the domain of $s$ contains two copies of each $2$-face $F$ in 
the target boundary component of $C_\phi$. They correspond to the source and target spaces $V_{F_s}$, $V_{F_t}$ of $Q_N^{r(F)}$. 
\smallskip

Finally, we provide an operator theoretic definition of the basic matrix dilogarithms (it is a straightforward rewriting of formula (32) in \cite{GT}), as well as their entries and those of the endomorphisms $Q_N$. Since it is the transposed tensor $\Hh_N^{red}(T_{C_\phi},\tilde b, \bw)^T$ mentioned in the Introduction (see formula \eqref {transposed}) that occurs in Theorem \ref{MAINTEO}, we consider the transposed endomorphisms of $Q_N$ and of the basic matrix dilogarithm instead:
$$Q_N^{T} : V_{F_s}' \to V_{F_t}'\ ,\ \Ll_N^T(\Delta,b,\bw): V_3'  \otimes V_1' \to V_2' \otimes V_0'$$
where $V_{F_s}' = V_{F_t}' = V_j'=(\mc^N)'$, the dual space of $\mc^N$, for every $j$. The matrix elements will be given with respect to the canonical basis of $\left((\C^N)'\right)^{\otimes 2}$.

\begin{remark}\label{identrem}{\rm If $\{e_j\}$ is the canonical basis of $\C^N$, and $\{e^j\}$ the dual basis of $(\C^N)'$, the identification map $\iota: \C^N \to (\C^N)'$, $\iota(e_j)=e^j$, extends to a canonical identification, also denoted by $\iota$, between $(\C^N)^{\otimes 2}$ and $\left( (\C^N)'\right)^{\otimes 2}$. Below we will rather deal with the endomorphisms
$$\iota^{-1}\circ Q_N^{T} \circ \iota \in {\rm End}\left(\C^N\right) \ ,\ \iota^{-1}\circ  \Ll_N^T(\Delta,b,\bw) \circ \iota \in {\rm End}\left((\C^N)^{\otimes 2}\right). $$
Obviously, they have the same matrix elements as $Q_N^{T}$ and $\Ll_N^T(\Delta,b,\bw)$ with respect to the canonical basis of $\C^N$ and $(\C^N)'$, and $(\C^N)^{\otimes 2}$ and $\left((\C^N)'\right)^{\otimes 2}$, respectively. Similarly we will consider
$$ \iota^{-1}\circ \Hh_N^{red}(T_{C_\phi},\tilde b, \bw)^T\circ \iota \in {\rm End}\left((\C^N)^{\otimes 2m}\right). $$
{\it For simplicity we will systematically omit the maps $\iota$ from the notations}.}
\end{remark}

The endomorphism $Q_N^T : V_{F_s} \to V_{F_t}$ is defined in the standard basis $\{e_j\}$ of (the copies $V_{F_s}$, $V_{F_t}$ of) $\mc^N$ by 
\begin{equation}\label{defQNT}
Q_N^T(e_k) = \frac{1}{\sqrt{N}} \sum_{j=0}^{N-1}q^{2kj+j^2}e_j.
\end{equation} 
Define the endomorphisms $A_0$, $A_1$, $A_2 \in {\rm End}(\mc^N)$ by 
\begin{equation}\label{standardmatrix}
A_0e_k=q^{2k}e_k, \ A_1e_k=e_{k+1}, \ A_2e_k=q^{1-2k}e_{k-1}.
\end{equation}
Note that $A_0A_1A_2 = q{\rm Id}_{\mc^N}$.
We stress that these operators arise in the representation of the {\it triangle algebra} (see Section \ref{standard} below);
their occurrence in the definition of the basic matrix dilogarithm is a key point to perform the computations at the end of the paper that
will establish the bridge between the QHI and the QT intertwiners.

Put 
$$g(x)=\prod_{j=1}^{N-1} (1-xq^{2j})^{j/N}, h(x)=\frac{g(x)}{g(1)}$$
$$\Upsilon = \frac{1}{N} \sum_{i,j=0}^{N-1} q^{-2ij-j} A_0^i \otimes (A_1A_0)^{-j}.$$
For any triple of quantum shape parameters $\bw=(\bw_0,\bw_1,\bw_2)$ and any endomorphism $U$ such that $U^N=-{\rm Id}$, set
$$\Psi_\bw(U) =  h(\bw_0) \sum_{i=0}^{N-1} U^i \prod_{j=1}^i \frac{-q^{-1}\bw_2}{1-\bw_0^{-1}q^{2j}}$$
where by convention we set the product equal to $1$ when $i=0$.  Then
\begin{equation}\label{Stensor}
\Ll_N^T(\Delta,b,\bw) = \Psi_\bw(-A_0A_1 \otimes A_1^{-1})\circ \Upsilon.
\end{equation}

Later we will need the following elementary but crucial fact, which belongs (in a different form) to Faddeev-Kashaev \cite{FK}. 
\begin{lem} Let $\bw_0$, $\bw_2$ be complex numbers such that $\bw_2^N = 1-\bw_0^{-N}$, and $U$ an endomorphism such that $U^N=-{\rm Id}$. Then $\Psi_\bw(U)$ defined as above is the unique endomorphism up to scalar multiplication which is a solution of the functional relation
\begin{equation}\label{cycdil2}
\Psi_\bw(q^2U) = (\bw_0+q^{-1}\bw_0\bw_2U)\Psi_\bw(U).
\end{equation}
\end{lem}
\Dim The hypothesis on $\bw_0$, $\bw_2$ implies that the summands of $\Psi_\bw(U)$ are periodic in $i$, with period $N$. Hence we get the same sum if $i$ ranges from $1$ to $N$. With this observation, we have
\begin{align*}
q^{-1}\bw_0\bw_2U\Psi_\bw(U) & = h(\bw_0) \sum_{i=0}^{N-1} U^{i+1} q^{-1}\bw_0\bw_2\prod_{j=1}^i \frac{-q^{-1}\bw_2}{1-\bw_0^{-1}q^{2j}} \\
 & = h(\bw_0) \sum_{i=0}^{N-1} U^{i+1} q^{-1}\bw_0\bw_2  \left(\frac{1-\bw_0^{-1}q^{2i+2}}{-q^{-1}\bw_2} \right)  \prod_{j=1}^{i+1} \frac{-q^{-1}\bw_2}{1-\bw_0^{-1}q^{2j}} \\
& =  h(\bw_0) \left(  - \bw_0 \sum_{i=0}^{N-1} U^{i+1} \prod_{j=1}^{i+1} \frac{-q^{-1}\bw_2}{1-\bw_0^{-1}q^{2j}}  + \sum_{i=0}^{N-1} (q^2U)^{i+1} \prod_{j=1}^{i+1} \frac{-q^{-1}\bw_2}{1-\bw_0^{-1}q^{2j}} \right)\\
& = - \bw_0 \Psi_\bw(U) + \Psi_\bw(q^2U).
\end{align*}
Hence $\Psi_\bw(U)$ solves the equation \eqref{cycdil2}. Conversely, $U^N=-{\rm Id}$ implies that $U$ is diagonalizable with eigenvalues in the set $\{-q^{2i},i=0,\ldots,N-1\}$. Hence $\Psi_\bw(U)$ is a polynomial in $U$ completely determined by its eigenvalues, which are of the form $\Psi_\bw(-q^{2i})$. By \eqref{cycdil2}  they are given for $i\geq 1$ by $\textstyle \Psi_\bw(-q^{2i}) = \prod_{j=1}^{i} (\bw_0-q^{2j-3}\bw_0\bw_2)\Psi_\bw(-1)$, and hence uniquely determined up to the choice of $\Psi_\bw(-1)$. Therefore \eqref{cycdil2} has a unique solution $\Psi_\bw(U)$ up to scalar multiplication. \hfill $\Box$

\begin{remark}{\rm The normalization factor $h(\bw_0)$ is chosen so that det$\left(\Psi_\bw(-A_0A_1 \otimes A_1^{-1})\right)=1$.}
\end{remark}

In particular, when $(\bw_0,\bw_1,\bw_2)$ is a triple of quantum shape parameters we have
\begin{equation}\label{cycdil}
\Psi_\bw(q^2U) = (\bw_0-\bw_1^{-1}U)\Psi_\bw(U).
\end{equation}
Finally the entries of $\Ll_N^T(\Delta,b,\bw)$ in the basis $\{e_k \otimes e_l\}_{k,l}$ of $\left(\mc^N\right)^{ \otimes 2}$ are as follows. Identify $\Z/N\Z$ with $\{0,1,\dots, N-1\}$. For every $a\in \Z$ set $\delta(a)=1$ if $a\equiv 0$ mod$(N)$, and $\delta(a)=0$ otherwise. For every $n\in \Z/N\Z$, consider the function $\omega(x,y|n)$ defined on the curve $\{(x,y\in \mc^2\vert \ x^N+y^N=1\}$ by
$$ \omega(x,y|0)=1, \omega(x,y|n)=\prod_{j=1}^n \frac{y}{1-xq^{-2j}}, \ n\geq
1.$$ 
Then, for every $i,j,k,l\in \Z/N\Z$ we have
$$\Ll_N^T(\Delta, b, \bw)_{i,j}^{k,l}= h(\bw_0) q^{-2kj-k^2}\omega(\bw_0, \bw_1^{-1}| i-k)\delta(i+j-l).$$ 

\subsection{Representations of the quantum Teichm\"uller spaces}\label{QREP} Unless stated differently, the results recalled in this
section are proved in \cite{B-L, B-B-L} or \cite{Filippo}. 
\medskip

Let $\lambda$ be an ideal triangulation of $S$, and $q$ as above. Recall that $m=-\chi(S)$. Put $n:=3m$ and fix an ordering
$e_1,\ldots,e_{n}$ of the edges $\lambda$. For all distinct $i$, $j$ set $\sigma_{ij} := a_{ij}-a_{ji} \in\{0,\pm 1,\pm 2\}$, where $a_{ij}$ is the number of times $e_i$ is on the {\it right} of $e_j$ in a triangle of $\lambda$, using the orientation of $S$. The {\it Chekhov-Fock algebra} $\TG_\lambda^q$ is the algebra over $\mc$ with generators $X_i^{\pm 1}$ associated to the edges $e_i$, and relations
\begin{equation}\label{or1}
X_iX_j=q^{2\sigma_{ij}}X_jX_i .
\end{equation} When $q=1$, $\TG_\lambda^1$ is just the algebra of Laurent polynomials $\mc[X_1^{\pm 1},\ldots,X_n^{\pm 1}]$, which is the ring of functions on the {\it classical} (enhanced) Teichm\"uller space generated by the exponential shear-bend coordinates on $\lambda$. 

\begin{figure}[ht]
\begin{center}
\includegraphics[width=7cm]{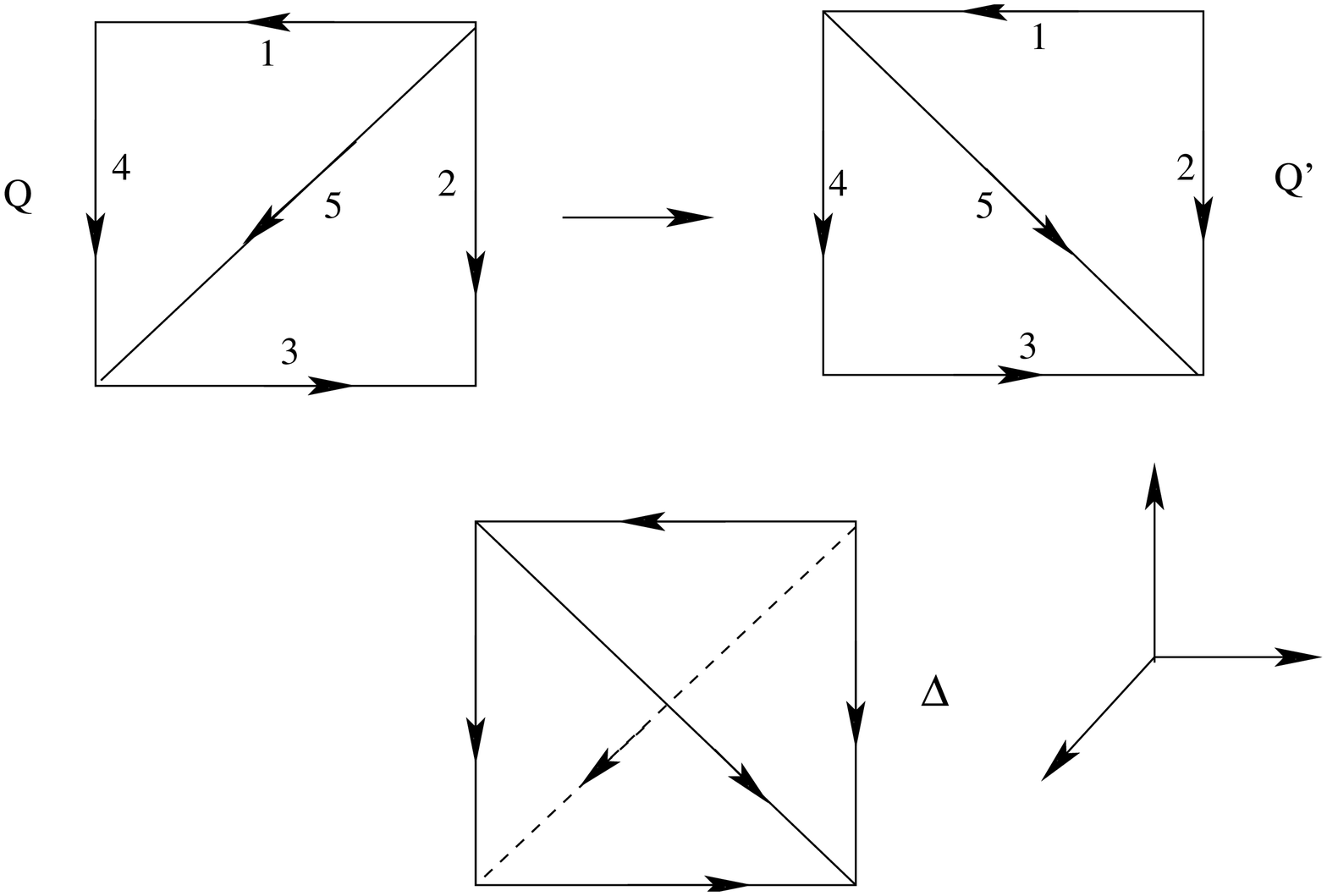}
\caption{\label{typefig2} }
\end{center}
\end{figure}

The algebra $\TG_\lambda^q$ has a well-defined fraction algebra $\hat \TG_\lambda^q$, and any diagonal exchange $\lambda \to \lambda'$ induces an isomorphism of algebras
$$\varphi_{\lambda\lambda'}^q: \hat \TG_{\lambda'}^q \rightarrow \hat \TG_\lambda^q.$$
The above definition of $\TG_\lambda^q$ works as well for any ideally triangulated punctured compact oriented surface $S$, possibly with boundary, where each boundary component is a union of edges.  In particular, take $S$ the ideal triangulated squares $Q$, $Q'$ in Figure \ref{typefig}. Number the edges of their triangulations $\lambda$, $\lambda'$ as in Figure \ref{typefig2}. In such a situation $\varphi_{\lambda\lambda'}^q$ has the form 
\begin{equation}\label{relq}
\begin{array}{ll}
\varphi_{\lambda\lambda'}^q(X_5') = X_5^{-1}\\
\varphi_{\lambda\lambda'}^q(X_1') = (1+qX_5)X_1& ,\ \varphi_{\lambda\lambda'}^q(X_3') = (1+qX_5)X_3\\ 
\varphi_{\lambda\lambda'}^q(X_2') = (1+qX_5^{-1})^{-1}X_2&  ,\ \varphi_{\lambda\lambda'}^q(X_4') = (1+qX_5^{-1})^{-1}X_4.
\end{array}
\end{equation}
It is proved in \cite{Bai0} (see also Theorem 1.22 in \cite{Filippo}) that this case of the ideal square determines, for {\it any} punctured compact oriented surface $S$, a unique family $\{\varphi_{\lambda\lambda'}^q\}_{\lambda\lambda'}$ of algebra isomorphisms $\varphi_{\lambda\lambda'}^q: \hat \TG_{\lambda'}^q \rightarrow \hat \TG_\lambda^q$ defined for all ideal triangulations $\lambda$ and $\lambda'$ of $S$, if the family satisfies certain natural properties with respect to composition, diffeomorphisms of $S$, and decomposition into ideally triangulated subsurfaces. 

Given any ideal triangulations $\lambda$, $\lambda'$ of $S$ and finite dimensional representation $\rho_\lambda: \TG_\lambda^q \rightarrow {\rm End}(V_\lambda)$, one says that $\rho_\lambda\circ \varphi_{\lambda\lambda'}^q$ {\it makes sense} if for every generator $X_i'\in \TG_{\lambda'}^q$ the element $\varphi_{\lambda\lambda'}^q(X_i')$ can be written as $P_iQ_i^{-1} \in \hat \TG_{\lambda}^q$, where $P_i$, $Q_i \in \TG_{\lambda}^q$ and $\rho_\lambda(P_i)$, $\rho_\lambda(Q_i)$ are invertible endomorphisms of $V_\lambda$. In such a case, $\rho_\lambda\circ \varphi_{\lambda\lambda'}^q(X_i') = \rho_\lambda(P_i)\rho_\lambda(Q_i)^{-1}$ is well-defined (ie. independent of the choice of the pair $(P_i,Q_i)$). 
\smallskip

By definition, the {\it quantum Teichm\"uller space} of $S$ is the quotient set $\textstyle \TG_S^q := (\coprod_{\lambda}\hat \TG_\lambda^q)/\sim$, where $\sim$ identifies the algebras $\hat \TG_\lambda^q$ and $\hat \TG_{\lambda'}^q$ by $\varphi_{\lambda\lambda'}^q$, for all ideal triangulations $\lambda$, $\lambda'$ of $S$. A {\it representation of $\TG_S^q$} is a family of representations
$$\rho:=\{\rho_\lambda: \TG_\lambda^q \rightarrow {\rm End}(V_\lambda)\}_{\lambda}$$
indexed by the set of ideal triangulations of $S$, such that $\rho_\lambda\circ \varphi_{\lambda\lambda'}^q$ makes sense for every $\lambda$, $\lambda'$ and is isomorphic to $\rho_{\lambda'}$. In fact, it is enough to check the isomorphisms $\rho_{\lambda'}\cong \rho_\lambda\circ \varphi_{\lambda\lambda'}^q$ whenever $\lambda$ and $\lambda'$ differ by a flip.

\subsubsection{Local  representations.}\label{localsection} The irreducible representations of $\TG_S^q$ were classified up to isomorphism in \cite{B-L}. The {\it local} representations of $\TG_S^q$ are special {\it reducible} representations defined as follows (see \cite{B-B-L,Filippo}). Define the {\it triangle algebra} $\Tt$ as the algebra over $\C$ with generators $Y_0^{\pm 1}$, $Y_1^{\pm 1}$, $Y_2^{\pm 1}$ and relations
\begin{equation}\label{or2}
Y_0Y_1=q^2Y_1Y_0,\ Y_1Y_2=q^2Y_2Y_1,\ Y_2Y_0=q^2Y_0Y_2.
\end{equation}
Fix an ordering $\tau_1,\ldots,\tau_{2m}$ of the triangles of $\lambda$. Order the abstract edges $e_{0}^j,e_{1}^j,e_{2}^j$ of $\tau_j$ so that the induced cyclic ordering is counter-clockwise with respect to the orientation of $S$. Associate to $\tau_j$ a copy $\Tt_j$ of the algebra $\Tt$, with generators denoted by $Y_0^j$,$Y_1^j$,$Y_2^j$, so that $Y_i^j$ is associated to the edge $e_{i}^j$. There is an algebra embedding 
\begin{equation}\label{embedi}
\mathfrak{i}_\lambda : \TG_\lambda^q\rightarrow  \Tt_1 \otimes \ldots \otimes \Tt_{2m}
\end{equation}
defined on generators by (we denote a monomial $\textstyle \otimes_j A_j$ in $\otimes_j\Tt_j$ by omitting  the terms $A_j=1$):
\begin{itemize}
\item If $e_i$ is an edge of two distinct triangles $\tau_{l_i}$ and $\tau_{r_i}$, and $e_{a_i}^{l_i}$, $e_{b_i}^{r_i}$ are the edges of $\tau_{l_i}$, $\tau_{r_i}$ respectively identified to $e_i$, then $\mathfrak{i}_\lambda(X_i) := Y_{a_i}^{l_i}\otimes Y_{b_i}^{r_i}$.
\item If $e_i$ is an edge of a single triangle $\tau_{k_i}$, and $e_{a_i}^{k_i}$, $e_{b_i}^{k_i}$ are the edges of $\tau_{k_i}$ identified to $e_i$, with $e_{a_i}^{k_i}$ on the right of $e_{b_i}^{k_i}$, then $\mathfrak{i}_\lambda(X_i) := q^{-1}Y_{a_i}^{k_i}Y_{b_i}^{k_i} = qY_{b_i}^{k_i}Y_{a_i}^{k_i}$.
\end{itemize}
A representation $\rho_\lambda$ of $\TG_\lambda^q$ is {\it local} if 
$$\rho_\lambda = (\rho_1\otimes \ldots \otimes \rho_{2m})\circ \mathfrak{i}_\lambda \colon \TG_\lambda^q \rightarrow {\rm End}(V_1\otimes \ldots \otimes V_{2m})$$ for some irreducible representations $\rho_j:\Tt \rightarrow {\rm End}(V_j)$ of $\Tt$, $j\in \{1,\ldots,2m\}$. So, a local representation is an equivalence class of tuples $(\rho_1,\ldots,\rho_{2m})$, where two tuples are equivalent if their restrictions to the subalgebra $\mathfrak{i}_\lambda(\TG_\lambda^q)$ of $\Tt_1 \otimes \ldots \otimes \Tt_{2m}$ define the same representation. Two local representations $\rho_\lambda=(\rho_1\otimes \ldots \otimes \rho_{2m})\circ \mathfrak{i}_\lambda$ and $\rho_\lambda'=(\rho_1'\otimes \ldots \otimes \rho_{2m}')\circ \mathfrak{i}_\lambda$ of $\TG_\lambda^q$ are isomorphic if there are linear isomorphisms $L_j\colon V_j\to V_j'$ such that for every $j=1,\ldots,2m$ and $Y \in \Tt_j$ we have
$$L_j \circ \rho_j(Y) \circ (L_j)^{-1} = \rho_j'(Y).$$
It is straightforward to check that this definition is independent of the choice of tuples $(\rho_1,\ldots,\rho_{2m})$, $(\rho_1',\ldots,\rho_{2m}')$. By definition a local representation
$\rho$ of $\TG_S^q$  is a representation formed by local representations $\rho_\lambda$.
\smallskip

\subsubsection{Isomorphism classes  and standard local representations}\label{standard} 
The isomorphism classes of irreducible representations of the triangle algebra $\Tt$ are parametrized by tuples of non zero scalars $(y_0,y_1,y_2,h)\in (\C^\bullet)^4$ such that $h^N=y_0y_1y_2$. The parameter $h$ is called the {\it load} of the class. The isomorphism class with parameters $(y_0,y_1,y_2,h)$ can be represented by {\it standard} representations $\rho:\Tt \to {\rm End}(\C^N)$, which have the form
\begin{equation}\label{paramloc}
\rho(Y_i)=\by_iA_i, \ i\in \{0,1,2\}
\end{equation}
where $\by_0,\by_1,\by_2\in \mc^*$ satisfy
$$\by_i^N=y_i,\ i\in \{0,1,2\},\ h=\by_0\by_1\by_2$$
and $A_0,A_1,A_2$ are the endomorphisms of $\mc^N$ defined in \eqref{standardmatrix}. 
\smallskip

Any local representation of $\TG_\lambda^q$ has dimension $N^{2m}$. The isomorphism class of a 
local representation $\rho_\lambda$ of $\TG_\lambda^q$ is determined by:
\begin{itemize}
\item a non zero complex {\it weight} $x_i$ associated to each edge of $\lambda$;
\item a $N$-th root $h$ of $x_1\ldots x_n$, called the {\it load}.
\end{itemize}
The weights $x_i$ and the load $h$ are such that
\begin{equation}\label{relparamloc}
\rho_\lambda(X_i^N) = x_i {\rm Id}_{V_1\otimes \ldots \otimes V_{2m}},\ \rho_\lambda(H) = h {\rm Id}_{V_1\otimes \ldots \otimes V_{2m}}
\end{equation}
where $H$ is a central element of $\TG_\lambda^q$, called {\it the principal central element}, given by \begin{equation}\label{punchdef}
H = q^{-\sum_{l<l'}\sigma_{ll'}} X_1\ldots X_n.
\end{equation}

Note that $h^N=x_1\ldots x_n$. It is straightforward to check that two local representations are isomorphic 
if and only if they are isomorphic as local representations. 
We call $(\{x_1,\ldots,x_n\},h)$ the {\it parameters} of $\rho_\lambda$. We say that a local representation 
$\rho_\lambda=(\rho_1\otimes \ldots \otimes \rho_{2m})\circ \mathfrak{i}_\lambda$ is standard if every $\rho_i$ is.
This notion naturally extends to representations of $\TG_S^q$.  
Every point of $(\mc^*)^{n}$ can be realized as the $n$-tuple of parameters $x_i$, $i=1,\dots,n$, of a standard local representation of $\TG_\lambda^q$. There is a one-to-one correspondence between isomorphism classes of local representations 
$\{\rho_\lambda: \TG_\lambda^q \rightarrow {\rm End}(V_1\otimes \ldots \otimes V_{2m})\}_{\lambda}$ of $\TG_S^q$ 
and families of parameters $\{(\{x_1,\ldots,x_n\}_\lambda,h_\lambda)\}_{\lambda}$ such that
\begin{equation}\label{coordch}
h_\lambda=h_{\lambda'}\ {\rm and}\ x_i=\varphi_{\lambda\lambda'}^1(x_i')
\end{equation}
for every $i=1,\ldots,n$ and any two ideal triangulations $\lambda$, $\lambda'$ of $S$ with edge weights $\{x_1,\ldots,x_n\}$, $\{x_1',\ldots, x_n'\}$ and loads $h_\lambda$, $h_{\lambda'}$ respectively. The edge weights $x_i$ of a local representation $\rho_\lambda: \TG_\lambda^q \rightarrow {\rm End}(V_\lambda)$ define its ``classical shadow'' $$sh(\rho_\lambda):\TG_{\lambda}^1 \rightarrow {\rm End}(\mc)$$ by taking the restriction of $\rho_\lambda$ to the subalgebra $\mc[X_1^{\pm N},\ldots, X_n^{\pm N}]\cong \TG_{\lambda}^1$ of the center of $\TG_\lambda^q$. This notion extends immediately to local representations $\rho=\{\rho_\lambda: \TG_\lambda^q \rightarrow {\rm End}(V_\lambda)\}_{\lambda}$ of $\TG_S^q$; the classical shadow
$$sh(\rho)=\{sh(\rho_\lambda): \TG_{\lambda}^1 \rightarrow {\rm End}(\mc)\}_{\lambda}$$ is a representation of the coordinate ring $\TG_S^1$ of the Teichm\"uller space (considered as a rational manifold with transition functions the maps $\varphi^1_{\lambda\lambda'}$). Every representation of $\TG_S^1$ is the shadow of $N$ local representations of $\TG_S^q$. The shadows of local representations encode the Zariski open subset of $X(S)$ made of the so-called {\it peripherically generic characters}. These include all the augmented characters of injective representations. 

The load $h$ has the following geometric interpretation. Let $[(\rr, \{\xi_{\Gamma}\}_{\Gamma \in \Pi})]$ be the augmented character of $S$ determined by $sh(\rho)$. Each 
$\Gamma \in \Pi$ is a group generated by the class in $\pi_1(S)$ of a small loop $m_j$ in $S$ going once and counter-clockwise around the $j$-th puncture, for some $j\in \{1,\ldots, r\}$ and some choice of basepoints. Since $\pi_1(S)$ is a free group,  $\rr\colon \pi_1(S) \rightarrow PSL(2,\mc)$ can be lifted to a homomorphism $\hat \rr\colon \pi_1(S) \rightarrow SL(2,\mc)$; the fixed point $\xi_\Gamma \in \mathbb{P}^1$ is then contained in an eigenspace of $\hat \rr(m_j)$, corresponding to an eigenvalue $a_j\in \mc$. Then
\begin{equation}\label{geomintpar}
h^N = (-1)^r a_1^{-1}\ldots a_r^{-1}.
\end{equation}
Finally, the decomposition into irreducible summands of a local representation $\rho_\lambda$ of $\TG_\lambda^q$ is (\cite{Tou})
\begin{equation}\label{decompirr}
\rho_\lambda = \oplus_{\rho_\lambda(\mu) \subset \rho_\lambda} \ \rho_\lambda(\mu)
\end{equation}
where the sum ranges over the set of all spaces $\rho_\lambda(\mu)$ formed by intersecting one eigenspace for each of the so-called (central) {\it puncture elements} of $\TG_\lambda^q$. Each space $\rho_\lambda(\mu)$ is also the {\it isotypical} component of an irreducible representation $\mu$ of $\TG_\lambda^q$, that is, the direct sum of all irreducible summands of $\rho_\lambda$ isomorphic to $\mu$. Every irreducible representation of $\TG_\lambda^q$ has an isotypical summand appearing in some local representation, and has multiplicity $N^g$ in it.

\subsection{Intertwiners of local representations}\label{Filippowork}
Let $(S,\lambda)$ and $q$ be as before. Any surface $R$ obtained by splitting $S$ along some (maybe all) edges of $\lambda$ inherits an orientation from $S$, and an ideal triangulation $\mu$ from $\lambda$. By gluing along edges backwards, one says that $(S,\lambda)$ is obtained {\it by fusion} from $(R,\mu)$. In such a case, every local representation $\eta_\mu =(\eta_1\otimes \ldots \otimes \eta_{2m})\circ \mathfrak{i}_\mu$ of $\TG_\mu^q$ determines a local representation $\rho_\lambda$ of $\TG_\lambda^q$ by setting $\rho_\lambda= (\eta_1\otimes \ldots \otimes \eta_{2m})\circ \mathfrak{i}_\lambda$. One says that $\eta_\mu$ {\it represents} $\rho_\lambda$. 

Given local representations $\eta = \{\eta_\mu \colon \TG_\mu^q \to {\rm End}(W_\mu)\}_\mu$,  $\rho = \{\rho_\lambda \colon \TG_\lambda^q \to {\rm End}(V_\lambda)\}_\lambda$ of $\TG_R^q$ and $\TG_S^q$ respectively, one says that {\it $\rho$ is obtained by fusion from $\eta$} if $\eta_\mu$ represents $\rho_\lambda$ for all ideal triangulations $\mu$ of $R$,  where $\lambda$ is the ideal triangulation of $S$ obtained by fusion from $\mu$.
\smallskip

In \cite{Filippo}, the following result is proved. The intertwiners $L_{\lambda\lambda'}^{\rho\rho'} \in \mathcal{L}_{\lambda\lambda'}^{\rho\rho'}$ in the statement are called {\it QT intertwiners}.
\begin{teo}\label{Filippoteo} There exists a collection $\{(\mathcal{L}_{\lambda\lambda'}^{\rho\rho'},\psi_{\lambda\lambda'}^{\rho\rho'})\}$, indexed by the couples of isomorphic local representations $\rho = \{\rho_\lambda \colon \TG_\lambda^q \to {\rm End}(V_\lambda)\}_\lambda$, $\rho' = \{\rho_\lambda' \colon \TG_\lambda^q \to {\rm End}(V_\lambda')\}_\lambda$ of $\TG_S^q$ and by the couples of ideal triangulations $\lambda$, $\lambda'$ of $S$, such that:
\smallskip

\noindent (1) $\mathcal{L}_{\lambda\lambda'}^{\rho\rho'}$ is a set of projective classes of linear isomorphisms $L_{\lambda\lambda'}^{\rho\rho'} \colon V_\lambda \to V_{\lambda'}$ such that for every $X'\in \TG_{\lambda'}^q$ we have
$$L_{\lambda\lambda'}^{\rho\rho'} \circ (\rho_\lambda \circ  \varphi^q_{\lambda\lambda'})(X')\circ (L_{\lambda\lambda'}^{\rho\rho'})^{-1}  = \rho_{\lambda'}'(X').$$

\noindent (2) $\psi_{\lambda\lambda'}^{\rho\rho'}\colon H_1(S;\mz/N\mz) \times \mathcal{L}_{\lambda\lambda'}^{\rho\rho'}\to \mathcal{L}_{\lambda\lambda'}^{\rho\rho'}$ is a free transitive action.
\smallskip

\noindent (3) Let $R$ be a surface such that $S$ is obtained by fusion from $R$. Let $\eta = \{\eta_\mu \colon \TG_\mu^q \to {\rm End}(W_\mu)\}_\mu$, $\eta' = \{\eta_\mu' \colon \TG_\mu^q \to {\rm End}(W_\mu')\}_\mu$ be two local representations of $\TG_R^q$ such that $\rho$, $\rho'$ are obtained respectively by fusion from $\eta$, $\eta'$. Then, for every ideal triangulations $\mu$, $\mu'$ of $R$, if $\lambda$, $\lambda'$ are the corresponding  ideal triangulations of $S$, there exists an inclusion map $j\colon \mathcal{L}_{\mu\mu'}^{\eta\eta'}\to \mathcal{L}_{\lambda\lambda'}^{\rho\rho'}$ such that for every $L\in \mathcal{L}_{\mu\mu'}^{\eta\eta'}$ and every $c\in H_1(R;\mz/N\mz)$ the following holds:
$$(j\circ  \psi_{\mu\mu'}^{\eta\eta'})(c,L) = \psi_{\lambda\lambda'}^{\rho\rho'}(\pi_*(c),j(L))$$
where $\pi\colon R\to S$ is the projection map.
\smallskip

\noindent (4) For every isomorphic local representations $\rho$, $\rho'$, $\rho''$ of $\TG_S^q$ and for every ideal triangulations $\lambda$, $\lambda'$, $\lambda''$ of $S$ the composition map
$$\fleche{\mathcal{L}_{\lambda\lambda'}^{\rho\rho'} \times \mathcal{L}_{\lambda'\lambda''}^{\rho'\rho''}}{\mathcal{L}_{\lambda\lambda''}^{\rho\rho''}}{(L_1,L_2)}{L_2\circ L_1}$$
is well-defined, and for all $c$, $d\in H_1(S;\mz/N\mz)$ it satisfies
$$\psi_{\lambda'\lambda''}^{\rho'\rho''}(d,L_2)\circ \psi_{\lambda\lambda'}^{\rho\rho'}(c,L_1) = \psi_{\lambda\lambda''}^{\rho\rho''}(c+d,L_2\circ L_1).$$
\end{teo}
It is also proved in \cite{Filippo} that any collection of intertwiners satisfying a weak form of conditions (3) and (4) (not involving the actions $\psi_{\lambda\lambda'}^{\rho\rho'}$) contains the collection $\{\mathcal{L}_{\lambda\lambda'}^{\rho\rho'}\}$. So the latter is minimal with respect to these conditions. 

Note that property (3) describes the behaviour of the intertwiners in $\mathcal{L}_{\lambda\lambda'}^{\rho\rho'}$ under cut-and-paste of subsurfaces along edges of $\lambda$, $\lambda'$. In particular, it implies that they can be decomposed into elementary intertwiners as follows. Any two ideal triangulations $\lambda$, $\lambda'$ of $S$ can be connected by a sequence $\lambda=\lambda_0 \rightarrow \ldots \rightarrow \lambda_{k+1}=\lambda'$ consisting of $k$ diagonal exchanges followed by an edge reindexing  $\lambda_k \rightarrow \lambda_{k+1}$ (with same underlying triangulation). Then, any projective class of intertwiners $[L^{\rho\rho'}_{\lambda\lambda'}] \in \mathcal{L}_{\lambda\lambda'}^{\rho\rho'}$ can be decomposed as
\begin{equation}\label{decompositionint}
L^{\rho\rho'}_{\lambda\lambda'}= L^{\rho\rho'}_{\lambda'\lambda'}\circ L^{\rho\rho}_{\lambda_{k}\lambda_{k+1}} \circ \dots \circ L^{\rho\rho}_{\lambda_0\lambda_{1}}
\end{equation}
where $L^{\rho\rho}_{\lambda_{i}\lambda_{i+1}} \in \mathcal{L}_{\lambda_i\lambda_{i+1}}^{\rho\rho}$ intertwins $\rho_{\lambda_{i}}$ and $\rho_{\lambda_{i+1}} \circ (\varphi^q_{\lambda_i\lambda_{i+1}})^{-1}$, related by the $i$-th diagonal exchange, for every $i\in \{0,\ldots,k-1\}$, and $L^{\rho\rho'}_{\lambda'\lambda'}$ intertwins $\rho_{\lambda'}$ and $\rho_{\lambda'}'$ on the triangulation $\lambda'$. In general the intertwiner $L^{\rho\rho'}_{\lambda\lambda'}$ depends on the choice of sequence $\lambda=\lambda_0 \rightarrow \ldots \rightarrow \lambda_{k+1}=\lambda'$, but the set $\mathcal{L}_{\lambda\lambda'}^{\rho\rho'}$ and the action $\psi_{\lambda\lambda'}^{\rho\rho'}$ do not.

\section{Proofs}\label{PF}
We are ready to prove our main Theorems \ref{MAINTEO}, \ref{MAINTEO2}, and Corollary \ref{cor}. Let us reformulate them by using the background material recalled in the previous sections. 

Let $M_\phi$ be a fibred cusped hyperbolic $3$-manifold realized as the mapping torus
of a pseudo Anosov diffeomorphism $\phi$ of a punctured surface $S$. Put $m=-\chi(S)>0$. Let $C_\phi$ be the associated cylinder. Let $T$ be
a layered triangulation of $M_\phi$, and $T_{C_\phi}$ the induced layered triangulation of $C_\phi$, 
with source boundary component the ideally triangulated surface $(S,\lambda)$, and target boundary component $(S,\lambda')$, where $\lambda' = \phi(\lambda)$ (as in Section  \ref{top-comb}).  

For every odd $N\geq 3$ and every primitive $N$-th root of unity $q$,  let $(T,\tilde b, \bw)$
be a QH layered triangulation of $M_\phi$, where  $\bw$ is a system of quantum shape parameters over $w\in A$, where $A$ is the 
subvariety  of the gluing variety $G(T,\tilde{b})$ as in  Proposition  \ref {gluing-var}. Let $(T_{C_\phi},\tilde b, \bw)$ be the induced QH triangulation
of $C_\phi$, and $\Hh_N^{red}(T_{C_\phi},\tilde b, \bw)\in {\rm End}\left((\C^N)^{\otimes 2m}\right)$ the QH operator defined by means of the
QH state sum carried by $(T_{C_\phi},\tilde b, \bw)$. Associated to it we have the transposed operator  
$\Hh_N^{red}(T_{C_\phi},\tilde b, \bw)^T\in {\rm End}\left(((\C^N)')^{\otimes 2m}\right)$ and 
$\iota^{-1}\circ  \Hh_N^{red}(T_{C_\phi},\tilde b, \bw)^T\circ \iota \in  {\rm End}\left((\C^N)^{\otimes 2m}\right)$ as in Section \ref{quantum-hyp}, still denoted by $\Hh_N^{red}(T_{C_\phi},\tilde b, \bw)^T$ (see Remark \ref{identrem}).

Now we re-state and prove our first main theorem, Theorem \ref{MAINTEO}.
\smallskip

{\bf First Main Theorem.}  {\it  (1) Every layered QH-triangulation $(T_{C_\phi},\tilde b,\bw)$ determines
  representations $\rho_\lambda$ and $\rho_{\lambda'}$ of $\TG_\lambda^q$ and $\TG_{\lambda'}^q$ respectively, belonging to a
  local representation $\rho$  of $\TG_{S}^q$ and such that $V_\lambda= V_{\lambda'}=(\C^N)^{\otimes 2m}$, and  
  $\rho_{\lambda}$ is isomorphic to $\phi^*\rho_{\lambda'}$. Moreover, the operator
  $\Hh_N^{red}(T_{C_\phi},\tilde b,\bw)^T$, considered as an element of ${\rm Hom}(V_\lambda,V_{\lambda'})$, is a QT intertwiner which intertwins the representations $\rho_\lambda$ and $\rho_{\lambda'}$.
  \smallskip

(2) For any other choice of weak branching $\tilde b'$,
  $\Hh_N^{red}(T_{C_\phi},\tilde b',\bw)^T$ intertwins local
  representations canonically isomorphic to $\rho_\lambda$,
  $\rho_{\lambda'}$ respectively.}
  
  \medskip
  
\begin{figure}[ht]
\begin{center}
\includegraphics[width=7cm]{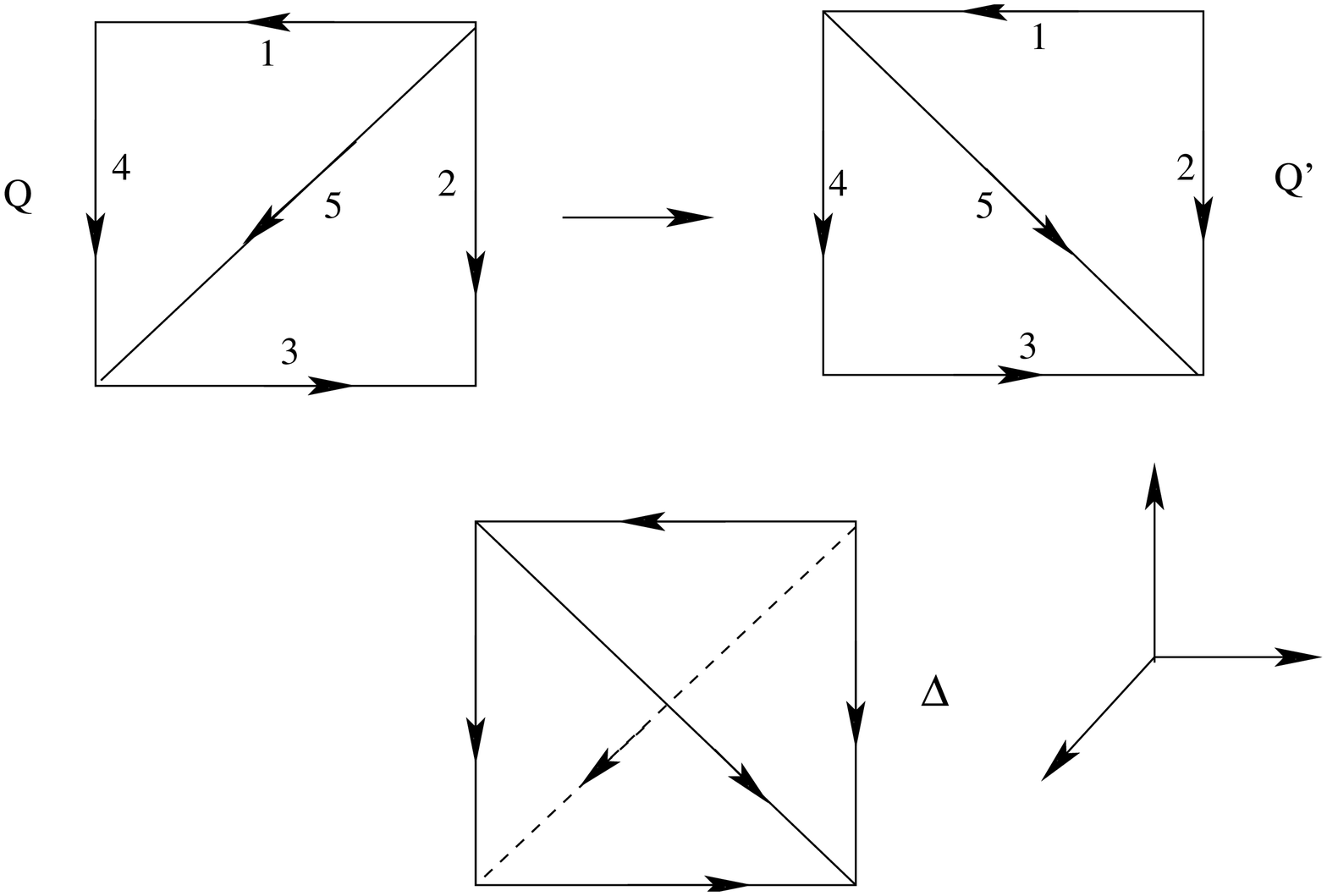}
\caption{\label{typefig3} }
\end{center}
\end{figure}

  \Dim We organize the proof in several steps:
\smallskip

\noindent {\bf Step 1.} Every transposed basic matrix dilogarithm $\Ll_N^T(\Delta,b,\bw)\in {\rm End}\left((\C^N)^{\otimes 2}\right)$ (see Section \ref{quantum-hyp})
 intertwins {\it standard} local representations of $\TG_\lambda^q$, $\TG_{\lambda'}^q$ (see Section  \ref{standard}), 
 where $\lambda$, $\lambda'$ are the ideal triangulations of the squares $Q$, $Q'$ in Figure \ref{typefig3}. 
 
\smallskip

\noindent {\bf Step 2.} Every transposed $2$-face operator $(Q_N^T)^{r(F)}: V_s \to V_t$ (see \eqref{defQNT}) intertwins representations of the triangle algebra 
$\rho_s\colon \Tt \to {\rm End}\left(\C^N\right)$, $\rho_t\colon \Tt \to {\rm End}\left(\C^N\right)$, associated to the abstract $2$-faces $F_s$, $F_t$ so that for each one the generator $Y_0$ labels the edge joining the lowest vertex to the middle one (with respect to the 
vertex ordering induced by the branching), and $Y_1$ labels the edge joining the lowest vertex to the biggest one. 
\smallskip

\noindent {\bf Step 3.} Every layered QH-triangulation $(T_{C_\phi},\tilde b,\bw)$ of the cylinder $C_\phi$, associated to a sequence
of diagonal exchanges $\lambda=\lambda_0\to \lambda_1 \to \dots \to \lambda_k=\lambda'$, determines standard local 
representations $\rho_j\colon\TG^q_{\lambda_j}\to {\rm End}\left( (\C^N)^{\otimes 2m}\right)$ associated to the surfaces $(S,\lambda_j)$, for every $j\in \{0,\ldots,k\}$. These representations belong to a local
representation $\rho$  
of the quantum Teichm\"uller space $\TG_S^q$. 

By Step 1, Step 2, and its actual definition by means of a QH state sum, $\Hh_N^{red}(T_{C_\phi},\tilde b,\bw)^T$ intertwins the local representations $\rho_0$ and $\rho_k$.
\medskip

Step (2) proves the claim (2) in the First Main Theorem; the three steps together prove (1). The details follow.

\subsubsection{Proof of Step 1}\label{localintert} Consider the squares $Q$, $Q'$ and the branched tetrahedron $(\Delta,b)$ in Figure \ref{typefig3}. Use the branching to order as $e_0$, $e_1$, $e_2$ the edges of any triangle of the squares $Q$, $Q'$, so that $e_0$ goes from the lowest vertex to the middle one, and $e_1$ goes from the lowest vertex to the biggest one. In such a situation there are two embeddings \eqref{embedi}, of the form $\mathfrak{i}_\lambda : \TG_\lambda^q\rightarrow  \Tt_3 \otimes \Tt_{1}$, $\mathfrak{i}_{\lambda'} : \TG_{\lambda'}^q\rightarrow  \Tt_2 \otimes \Tt_{0}$, where $\Tt_i$ is the copy of the triangle algebra associated to the $i$-th $2$-face (ie. opposite to the $i$-th vertex) of $(\Delta,b)$, and
$$\mathfrak{i}_\lambda(X_1) = Y_0^3 \otimes 1, \mathfrak{i}_\lambda(X_2) = 1 \otimes Y_1^1, \mathfrak{i}_\lambda(X_3) = 1\otimes Y_2^1, \mathfrak{i}_\lambda(X_4) = Y_2^3\otimes 1, \mathfrak{i}_\lambda(X_5) = Y_1^3\otimes Y_0^1.$$
$$\mathfrak{i}_{\lambda'}(X_1') =Y_0^2 \otimes 1,\mathfrak{i}_{\lambda'}(X_2') = Y_1^2\otimes 1,\mathfrak{i}_{\lambda'}(X_3') = 1\otimes Y_2^0,\mathfrak{i}_{\lambda'}(X_4') = 1 \otimes Y_0^0,\mathfrak{i}_{\lambda'}(X_5') = Y_2^2 \otimes Y_1^0.$$ 
Note that $Y_j^i$ is associated to the edge $e_j$ of the $i$-th $2$-face of $(\Delta,b)$, that we have specified above. Recall the notion of standard local representation in Section
\ref{standard}.

\begin{prop}\label{Sintertwins} Let $\rho_\lambda = (\rho_3\otimes \rho_{1})\circ \mathfrak{i}_\lambda \colon \TG_\lambda^q \rightarrow {\rm End}(V_3\otimes V_{1})$ be a standard local representation, and $\by_j^i$  be such that $\rho_i(Y_j^i) =\by_j^iA_j$, for $i\in \{1,3\}$, $j\in \{0,1,2\}$. Let $\bw$ be a triple of quantum shape parameters of the branched tetrahedron $(\Delta,b)$ such that $\bw_2 = -q \by_1^3\by_0^1$. Then, for all $X\in \TG_\lambda^q$ we have 
\begin{equation}\label{intprop}
\Ll_N^T(\Delta,b,\bw) \circ \rho_\lambda(X)  =  \rho_{\lambda'}((\varphi^q_{\lambda\lambda'})^{-1}(X))  \circ  \Ll_N^T(\Delta,b,\bw)
\end{equation}
where $\rho_{\lambda'} = (\rho_2\otimes \rho_{0})\circ \mathfrak{i}_\lambda \colon \TG_{\lambda'}^q \rightarrow {\rm End}(V_2\otimes V_{0})$ is the standard representation given by
\begin{equation}\label{relexdiag}
\by^2_0 = \bw_0^{-1} \by^3_0\quad , \quad \by^2_1 = \bw_1^{-1} \by^1_1\quad ,\quad \by^0_0 = \bw_1^{-1}\by^3_2\quad, \quad \by^0_2 = \bw_0^{-1} \by^1_2\quad , \quad \by^2_2\by^0_1 = -q\bw_2^{-1}.
\end{equation}
\end{prop}
\begin{remark}{\rm Note that $\bw_2$ and $w_1 = \bw_1^N = (1-\bw_2^N)^{-1}$ and $w_0= 1-w_1^{-1}$ are determined by $\rho_\lambda$, but the choice of $N$-th root $\bw_1$ (or $\bw_0$) is free; this choice determines $\bw$, and hence $\rho_{\lambda'}$, completely by the relation $\bw_0\bw_1\bw_2=-q$.  Also, we have $\rho_\lambda(X_5) = \by_1^3\by_0^1 A_1 \otimes A_0 = -q^{-1}\bw_2A_1 \otimes A_0$, $\rho_{\lambda'}(X_5') = -q\bw_2^{-1}A_2 \otimes A_1$, so that
$$\rho_\lambda(X_5^N) = -w_2 I_N \otimes I_N, \rho_{\lambda'}((X_5')^N)= -w_2^{-1} I_N \otimes I_N.$$ }
\end{remark}

\Dim It is enough to check \eqref{intprop} for $X\in \{X_1,\ldots,X_5\}$, using the relations \eqref{relq}. Let us do the cases $X=X_5$ and $X_4$, the other cases being respectively similar to this last. For $X=X_5$ recall that  $\varphi_{\lambda\lambda'}^q(X_5')=X_5^{-1}$; in this case  the identity \eqref{intprop} reads
$$\Ll_N^T(\Delta,b,\bw) \circ (\by^1_0\by^3_1 A_1 \otimes A_0)  =  
(\by^0_1\by^2_2 A_2 \otimes A_1)^{-1} \circ  \Ll_N^T(\Delta,b,\bw).$$
Consider the factorization formula \eqref{Stensor}. We have $(A_2 \otimes A_1)^{-1} = qA_1A_0 \otimes A_1^{-1}$, so it commutes with $\Psi_\bw(- A_0A_1 \otimes A_1^{-1})$ and we get
$$\begin{array}{l}
(\by^0_1\by^2_2 A_2 \otimes A_1)^{-1} \circ  \Ll_N^T(\Delta,b,\bw) = \\
\hspace*{0.5cm}  = \Psi_\bw(-A_0A_1 \otimes A_1^{-1})\circ  (q(\by^0_1\by^2_2)^{-1} A_1A_0 \otimes A_1^{-1}) \circ \left( \frac{1}{N} \sum_{i,j=0}^{N-1} q^{-2ij-j} A_0^i \otimes (A_1A_0)^{-j} \right)\\
\hspace*{0.5cm} = \Psi_\bw(-A_0A_1 \otimes A_1^{-1})\circ \left( \frac{1}{N} \sum_{i,j=0}^{N-1} q^{-2(i+1)(j+1)-(j+1)+3} A_0^i \otimes (A_1A_0)^{-j} \right) \circ \\
\hspace*{11cm} \circ (q(\by^0_1\by^2_2)^{-1} A_1A_0 \otimes A_1^{-1})\\
\hspace*{0.5cm} = \Ll_N^T(\Delta,b,\bw) \circ  (q^3 A_0^{-1}\otimes A_1A_0)\circ (q(\by^0_1\by^2_2)^{-1} A_1A_0 \otimes A_1^{-1})\\ \\
\hspace*{0.5cm} =  \Ll_N^T(\Delta,b,\bw) \circ  ((\by^0_1\by^2_2)^{-1} A_1 \otimes A_0).
\end{array}$$
Hence the identity \eqref{intprop} holds true for $X=X_5$ whenever
\begin{equation}\label{firstrel}
(\by^0_1\by^2_2)^{-1}  = \by^1_0\by^3_1 = -q^{-1}\bw_2.
\end{equation}
For $X=X_4$ recall that $\varphi_{\lambda\lambda'}^q(X_4')= (1+qX_5^{-1})^{-1}X_4$; in this case  the identity \eqref{intprop} reads
\begin{equation}\label{intprop4}
\Ll_N^T(\Delta,b,\bw) \circ (\by^3_2 A_2 \otimes I_N)  =  (1+q \by^2_2\by^0_1 A_2\otimes A_1) \circ (I_N \otimes \by^0_0 A_0) \circ  
 \Ll_N^T(\Delta,b,\bw).
\end{equation}
Now we have: 
\begin{align*}
\Ll_N^T(\Delta,b,\bw) \circ & (\by^3_2 A_2 \otimes I_N)  \\
&  = \Psi_\bw(-A_0A_1 \otimes A_1^{-1})\circ (\by^3_2 A_2 \otimes I_N)  \circ \left( \frac{1}{N} \sum_{i,j=0}^{N-1} q^{-2i(j+1)-j} A_0^i \otimes (A_1A_0)^{-j} \right)\\  & = \Psi_\bw(-A_0A_1 \otimes A_1^{-1})\circ (q\by^3_2 A_2 \otimes A_1A_0) \circ \Upsilon\\ & = (q\by^3_2 A_2 \otimes A_1A_0) \circ \Psi_\bw(q^2(-A_0A_1 \otimes A_1^{-1}))\circ \Upsilon \\ & = (q\by^3_2 A_2 \otimes A_1A_0) \circ (\bw_0+q^2\bw_1^{-1} A_1A_0 \otimes A_1^{-1}) \circ \Psi_\bw(-A_0A_1 \otimes A_1^{-1})\circ \Upsilon\\ & = (1+\bw_0\bw_1 (A_1A_0)^{-1} \otimes A_1) \circ  (\bw_1^{-1} A_1A_0 \otimes A_1^{-1}) \circ \\ & \hspace*{6cm} \circ  (q\by^3_2 A_2 \otimes A_1A_0) \circ \Ll_N^T(\Delta,b,\bw) 
 \\ & = (1-q^2\bw_2^{-1} A_2 \otimes A_1) \circ (I_N \otimes \bw_1^{-1}\by^3_2 A_0) \circ \Ll_N^T(\Delta,b,\bw).
\end{align*}
Note that we used the relation \eqref{cycdil} in the fourth equality. Hence, using \eqref{firstrel} we see that \eqref{intprop} holds true for $X=X_4$ whenever $\by^0_0 = \bw_1^{-1}\by^3_2$. The other cases $X=X_1$, $X_2$ or $X_3$ are similar. \hfill $\Box$

\subsubsection{Proof of Step 2} Consider QH-tetrahedra $(\Delta,b,\bw)$, $(\Delta',b',\bw')$ glued along a $2$-face $F$. Denote as usual by $F_s$, $F_t$ the abstract $2$-faces of $\Delta$, $\Delta'$ corresponding to $F$. Recall that the vertices $v_0$, $v_1$, $v_2$ of $F_s$ and $F_t$ are ordered by the branchings $b$ and $b'$ respectively; the gluing $F_s \to F_t$ is encoded by an even permutation $p:=\sigma^{r(F)}$ of the vertices, where $\sigma$ is the $2$-cycle $(v_0v_1v_2)$ and $r(F) \in \mz/3\mz$. 

Assume that a standard representation of the triangle algebra $\rho_s\colon \Tt \to {\rm End}(\C^N)$ is given on $F_s$, where the generators $Y_0$ labels the edge $[v_0,v_1]$, and $Y_1$ labels the edge $[v_0,v_2]$. Using the gluing, $\rho_s$ is the pull-back of a standard representation $\rho_t\colon \Tt \to {\rm End}(\C^N)$ on $F_t$, where now $Y_0$ labels the edge $[v_{p(0)},v_{p(1)}]$, and $Y_1$ labels the edge $[v_{p(0)},v_{p(2)}]$. So $Y_i$ on $F_s$ corresponds to $Y_{p^{-1}(i)}$ on $F_t$. 
\smallskip

We have:

\begin{lem} \label{changeb} The endomorphism $(Q_N^T)^{r(F)}$ intertwins $\rho_s$ and $\rho_t$. Namely, for all $X$ in $\Tt$ we have 
$(Q_N^T)^{r(F)}\circ \rho_s(X) = \rho_t(X) \circ (Q_N^T)^{r(F)}$.
\end{lem}
\Dim It is enough to check this on generators, where it follows from $Q_N^TA_i(Q_N^T)^{-1} = A_{i-1}$ (indices {\rm mod}$(3)$). This is straightforward to check.\cvd
\medskip

Note that $Q_N$ has order $3$ up to multiplication by $4$-th roots of $1$ (see \cite{AGT}, Lemma 7.3: 
if $N=2n+1$, we have $Q_N^3 = \phi_N^{-1}I_N$, where $\textstyle \phi_N = \left(\frac{n+1}{N}\right)$ if $N \equiv 1$ mod$(4)$, and $\textstyle \phi_N = \left(\frac{n+1}{N}\right) i$ if $N \equiv 3$ mod$(4)$).

\subsubsection{Proof of Step 3}  Let $(T_{C_\phi},\tilde b,\bw)$ be a layered QH-triangulation of the cylinder $C_\phi$ as usual. 
Recall that for every $j\in \{0,\dots, k\}$, the triangulated surface $(S,\lambda_j)$ inherits a weak branching $\tilde b_j$ from its positive side in $(T,\tilde b,\bw)$. Order as $\tau_1^j,\ldots,\tau_{2m}^j$ the branched abstract triangles of $(\lambda_j,\tilde b_j)$. For every $i\in \{1,\dots, 2m\}$, order as $e_0^i$, $e_1^i$, $e_2^i$ the edges of the triangle $\tau_i^j$, as described before Proposition \ref{Sintertwins} and Lemma \ref{changeb}. Label $e_k^i$ with the generator $Y_k^i$ of a copy of the triangle algebra $\Tt$. 
\smallskip

Assume that for some $j\in \{0,\dots, k\}$ we are given a standard local representation
$$\rho_j=(\rho_{1}^j\otimes\ldots \otimes \rho_{2m}^j) \circ i_{\lambda_j}\colon\TG^q_{\lambda_j}\to {\rm End}( (\C^N)^{\otimes 2m})$$
where the standard representations $\rho_{i}^j: \Tt\to {\rm End}(\C^N)$ are associated to the triangles $\tau_i^j$. As in \eqref{paramloc}, put $\rho_{i}^j(Y_k) :=\by_{i,k}^{j} A_k$. By the very definition of $\mathfrak{i}_{\lambda_j}$ (see \eqref{embedi}), for any edge $e_s$ of $\lambda_j$, with edge generator $X_{s}\in \TG^q_{\lambda_j}$, we have (as usual we denote a monomial $\textstyle \otimes_j A_j$ in $\Tt^{\otimes 2m}$ by omitting  the terms $A_j=1$):
\begin{itemize}
\item[(a)] $\rho_j(X_s) = \by_{l_i,a_i}^{j}\by_{r_i,b_i}^{j} A_{a_i}^{l_i} \otimes A_{b_i}^{r_i}$, if $e_s$ is an edge of two distinct triangles $\tau_{l_i}^j$ and $\tau_{r_i}^j$, and $e_{a_i}^{l_i}$, $e_{b_i}^{r_i}$ are the edges of $\tau_{l_i}^j$, $\tau_{r_i}^j$ respectively identified to $e_s$.
\item[(b)] $\rho_j(X_s) = q^{-1} \by_{k_i,a_i}^{j}\by_{k_i,b_i}^{j} A_{a_i}^{k_i}A_{b_i}^{k_i}$, if $e_s$ is an edge of a single triangle $\tau_{k_i}^j$, and $e_{a_i}^{k_i}$, $e_{b_i}^{k_i}$ are the edges of $\tau_{k_i}^j$ identified to $e_s$, with $e_{a_i}^{k_i}$ on the right of $e_{b_i}^{k_i}$.
\end{itemize}
Consider the diagonal exchange $\lambda_j\to \lambda_{j+1}$, and the corresponding triangulated (abstract) squares $Q$ and $Q'$, as in Figure \ref{typefig3}. The standard representations $\rho_i^j: \Tt\to {\rm End}(\C^N)$ associated to the triangles of $Q$ define a standard local representation $\rho_{Q}$ of $\TG^q_{\lambda_j|Q}$. Denote by $\rho_{Q'}$ the standard local representation of $\TG^q_{\lambda_{j+1}|Q'}$ defined from $\rho_{Q}$ by the relations \eqref{relexdiag}, where $(\Delta,b,\bw)$ is the QH-tetrahedron of $(T_{C_\phi},\tilde b,\bw)$ bounded by $Q$ and $Q'$; we describe below the relations between $\bw$ and $\rho_j$ (and hence $\rho_Q$). The representation $\rho_{Q'}$ extends (by fusion) to a standard local representation
$$\rho_{j+1}=(\rho_{1}^{j+1}\otimes\ldots \otimes \rho_{2m}^{j+1}) \circ i_{\lambda_j}\colon\TG^q_{\lambda_{j+1}}\to {\rm End}( (\C^N)^{\otimes 2m})$$
by stipulating that $\rho_{i}^{j+1} = \rho_{i}^{j}$ on the common triangles of $(S,\lambda_j)$ and $(S,\lambda_{j+1})$. Since $(S,\lambda_{j+1})$ is given the weak-branching $\tilde b_{j+1}$ induced by its positive side in $(T_{C_\phi},\tilde b,\bw)$, for every triangle $\tau_i^{j+1}$ of $(S,\lambda_{j+1})$ where $\tilde b_{j+1}$ differs from the weak-branching given by $(S,\lambda_j)$ or $Q'$, replace the corresponding representation $\rho_{i}^{j+1}$ by $(Q_N^T)^{r(F)} \circ \rho_i^{j+1}(X) \circ (Q_N^T)^{-r(F)}$, as in Lemma \ref{changeb}. For simplicity, we keep however the same notation for $\rho_{j+1}$. Define
$$L^{\rho_j\rho_{j+1}}_{\lambda_j\lambda_{j+1}}\colon {\rm End}((\mc^N)^{\otimes 2m}) \rightarrow {\rm End}((\mc^N)^{\otimes 2m})$$
as the operator acting as the identity on all factors except on the representation space $\mc^N\otimes \mc^N \hookrightarrow (\mc^N)^{\otimes 2m}$ of $\rho_{Q}$, where it acts by 
$\left((Q_N^T)^{r(F_2)} \otimes  (Q_N^T)^{r(F_0)}\right)\circ \Ll_N^T(\Delta,b,\bw)$. Here, as usual we denote by $F_2$ and $F_0$ the $2$-faces supporting the target space of $\Ll_N^T(\Delta,b,\bw)$.  By construction, for every $X\in \TG^q_{\lambda_{j+1}}$ we have
\begin{equation}\label{intertwj}
L^{\rho_j\rho_{j+1}}_{\lambda_j\lambda_{j+1}} \circ \rho_{j}\circ \varphi^q_{\lambda_j\lambda_{j+1}}(X')\circ (L^{\rho_j\rho_{j+1}}_{\lambda_j\lambda_{j+1}})^{-1} =  \rho_{j+1}(X).
\end{equation}
Finally, assuming that a standard local representation $\rho_j$ is given for $j=0$, by working as above we can define inductively a sequence of standard local
representation $\rho_j$, $j=0,\dots , k$, and set
$$L^{\rho_0,\rho_k}_{\lambda_0, \lambda_k}:= 
L^{\rho_{k-1}\rho_{k}}_{\lambda_{k-1}\lambda_{k}}\circ \dots \circ L^{\rho_{0}\rho_{1}}_{\lambda_{0}\lambda_{1}}.$$
Comparing with \eqref{Ssumformula} we see that
$$L^{\rho_0,\rho_k}_{\lambda_0, \lambda_k} = \Hh_N^{red}(T_{C_\phi},\tilde b, \bw)^T.$$

Now we have to specialize the choice of $\rho_0$, structurally related to the triangulation $(T_{C_\phi},\tilde b,\bw)$. Similarly to the ``classical" lateral shape 
parameters $W_j^+(e)$ (Section  \ref{class-hyp}), put:
\begin{defi}\label{qsbcoordinates} {\rm For every edge $e$ of $(S,\lambda_j)$, the {\it lateral quantum shape parameter} $\bW_j^+(e)$ is the product of the quantum shape parameters $\bw(E)$ of the abstract edges $E\rightarrow e$ carried by tetrahedra lying on the positive side of $(S,\lambda_j)$.  The {\it quantum shear-bend coordinates} of $(S,\lambda_j)$ are the scalars $$\bx^j(e):=-q^{-1}\bW_j^+(e).$$}
\end{defi} 

Denote by $\lambda_j^{(1)}$ the set of edges of the ideal triangulation $\lambda_j$. 
\smallskip

By Proposition \ref{shear-bend}, the quantum shear-bend coordinates $\bx^j(e)$, $e\in \lambda_j^{(1)}$, form a system of $N$-th roots of the shear-bend coordinates $x^j(e)$ of $\lambda_j$. By fixing $w$ and varying $\bw$ over $w$ in the tuple $(T_{C_\phi},\tilde b,\bw)$, one obtains a family of systems $\{\bx^j(e)\}_{e\in \lambda_j^{(1)}}$ such that $\bx^0(e)=\bx^k(\phi(e))$ for every edge $e$ of $\lambda=\lambda_0$. 

The QH-triangulation $(T_{C_\phi},\tilde b,\bw)$ determines $\bx^j(e_s)$ for all $j\in \{0,\ldots,k\}$ and every edge $e_s$ of $\lambda_j$. Consider the case $j=0$, and the regular map
$$\pi_0\colon \{\by^0 = (\by_{i,k}^{0})_{i,k}\} \to \{\bx^0 = (\bx^0(e_s))_{e_s\in \lambda_0^{(1)}}\}$$ defined by
\begin{equation}\label{relloc}
\bx^0(e_s) = \by_{l_i,a_i}^{0}\by_{r_i,b_i}^{0}\quad  {\rm or}\quad \by_{k_i,a_i}^{0}\by_{k_i,b_i}^{0}
\end{equation}
according to the cases (a) or (b) described at the beginning of Step 3. A simple dimensional count shows that $\pi_0$ is surjective, with generic fibre isomorphic to $(\mc^*)^{3m}$. Then, take a point $\by^0\in \pi^{-1}_0(\bx^0)$. We get
\begin{equation}\label{relloc2}
\rho_0(X_s) = \bx^0(e_s) A_{a_i}^{l_i} \otimes A_{b_i}^{r_i} \quad  {\rm or}\quad \rho_0(X_s) = q^{-1} \bx^0(e_s) A_{a_i}^{k_i}A_{b_i}^{k_i}
\end{equation}
again according to the cases (a) or (b) above. We have to check that the sequence of local representations $\rho_j$ constructed by starting from $\rho_0$ is consistent, that is, compatible with $\bw$ under the diagonal exchanges. Using the notations of Figure \ref{typefig3},  by the tetrahedral and edge relations satisfied by the quantum shape parameters, for all $j\in \{0,\ldots,k\}$ we have
\begin{equation}\label{relqsb}
\bx^j(e_5)=-q^{-1}\bw_2, \bx^{j+1}(e_5)=-q\bw_2^{-1}
\end{equation}
and
\begin{equation}\label{relqsb2}
\begin{array}{ccc}
\bx^{j+1}(e_1)=\bx^{j}(e_1)\bw_0^{-1}& , & \bx^{j+1}(e_2)=\bx^{j}(e_2)\bw_1^{-1}\\
\bx^{j+1}(e_3)=\bx^{j}(e_3)\bw_0^{-1}& , & \bx^{j+1}(e_4)=\bx^{j}(e_4)\bw_1^{-1}.
\end{array}
\end{equation}
In particular, \eqref{relqsb} for $j=0$ is consistent with the relation $\bw_2 = -q \by_1^3\by_0^1$ required by Proposition \ref{Sintertwins} (where $\by_1^3 = \by_{r_i,b_i}^{0}$ and $\by_0^1 = \by_{l_i,a_i}^{0}$). Then, comparing \eqref{relqsb2} and \eqref{relexdiag} we see that for every $j\in \{0,\ldots,k\}$, $\rho_j$ is indeed defined by identities like \eqref{relloc2}, with $\bx^0$ replaced by $\bx^j$.

\smallskip By \eqref{intertwj}, the representations $\rho_0,\ldots,\rho_k$ belong to a local representation $\rho=\{\rho_\lambda: \TG_\lambda^q \rightarrow {\rm End}(V_\lambda)\}_{\lambda}$ of $\TG^q_S$. By \eqref{relparamloc}, \eqref{relloc2} and the similar identities for $\rho_j$, we see that the parameters of $\rho_j$ are $$(\{x^j(e)\}_{e\in \lambda_j^{(1)}}, h)$$ where $x^j(e) = (\bx^j(e))^N$ and 
\begin{equation}\label{load2}
h =\prod_{e\in \lambda_j^{(1)}} \bx^j(e).
\end{equation}
For this identity, first note that the right hand side, say $h^j$, does not depend on $j$. Indeed, by \eqref{relqsb}, \eqref{relqsb2} and the tetrahedral and edge relations for $\bw_0$, $\bw_1$ and $\bw_2$, we see that $$\textstyle \prod_{i=1}^5 \bx^j(e_i) = \prod_{i=1}^5 \bx^{j+1}(e_i)$$
under a diagonal exchange. Since the quantum shear-bend coordinates of edges not touched by a diagonal exchange are not altered, we deduce that $h^j=h^{j+1}$ for every $j=0,\dots, k-1$. Now we can check \eqref{load2} as follows. Consider a triangle $\tau^j_s$ of $\lambda_j$, with edges $X_i,X_j,X_k$, where $i<j<k$. The weak branching $\tilde b$ provides a bijection $p\colon \{i,j,k\}\mapsto \{0,1,2\}$. Then $\mathfrak{i}_{\lambda_j}(H)$ is a tensor product of monomials associated to the triangles of $\lambda_j$, and the monomial for $\tau^j_s$ is $q^{-\sigma_{\vert \tau^j_s}}Y_{p(i)}^s Y_{p(j)}^s Y_{p(k)}^s$ where $-\sigma_{\vert \tau^j_s}$ is the contribution to the summand $-\sigma_{ij}-\sigma_{ik}-\sigma_{jk}$ of $\textstyle -\sum_{l<l'}\sigma_{ll'}$ coming from the triangle $\tau^j_s$. Noting that the product $Y_0^sY_1^sY_2^s$ is invariant under cyclic permutations, one checks immediately that $q^{-\sigma_{\vert \tau_j^l}}Y_{p(i)}^s Y_{p(j)}^s Y_{p(k)}^s = q^{-1} Y_0^sY_1^sY_2^s$. So  
$$\mathfrak{i}_{\lambda_j}(H) = \otimes_{s=1}^{2m} (q^{-1}Y_0^sY_1^sY_2^s)$$
and hence
$$\rho_j(H) = \otimes_{s=1}^{2m} (\by_0^s\by_1^s\by_2^sq^{-1}A_0^sA_1^sA_2^s) = \bx^j(e_1)\ldots \bx^j(e_n) {\rm Id}_{(\mc^N)^{\otimes 2m}}.$$

Summarizing our discussion, we have proved the following result. It concludes the proof of Step 3 of our First Main Theorem, and hence of Theorem \ref{MAINTEO} in the introduction.
\begin{prop}\label{QHIlocrep} 
 The QH triangulation $(T_{C_\phi},\tilde b,\bw)$ determines a sequence of standard local representations $\rho_j\colon\TG^q_{\lambda_j}\to {\rm End}( (\C^N)^{\otimes 2m})$ associated to the triangulated surfaces $(S,\lambda_j)$,  with parameters $$\textstyle (\{x^j(e)\}_{e\in \lambda_j^{(1)}},h),\ h=\prod_{e\in \lambda_j^{(1)}}\bx^j(e),$$ 
 and belonging to a local representation $\rho$ of  $\TG^q_{S}$. Moreover, $\Hh_N^{red}(T_{C_\phi},\tilde b,\bw)^T$ intertwins the representations $\rho_0$ and $\rho_k$.
\end{prop} 

Let us consider now Theorem \ref{MAINTEO2}. For the convenience of the reader we re-state it.
We keep the usual setting. Recall that $M_\phi$ is the interior of a compact manifold $\bar M_\phi$ with boundary made by tori.
\medskip

{\bf Second Main Theorem.}
  {\it There is a neighborhood of $i^*(\rG_h)$ in $i^*(X_0(M_\phi)) \subset
  X(S)$ such that:
  
  (1) For any isomorphism class of local representations of $\TG_S^q$
  whose holonomy lies in this neighborhood, there are: 
  \begin{itemize}
  \item A representative $\rho$ of the class, and representations $\rho_\lambda$, $\rho_{\lambda'}$ belonging to $\rho$ and acting on $(\C^N)^{\otimes 2m}$,
  \item a layered QH triangulation $(T,\tilde b, \bw)$ of $M_\phi$ such that $\Hh_N^{red}(T,\tilde b, \bw)=\Hh_N^{red}(M_\phi,\rG,\kappa)$, for some weight $\kappa\in H^1(\partial \bar M_\phi,\mc^*)$ and augmented character $\rG$ of $\pi_1(M_\phi)$ such that the restriction of $\rG$ to $\pi_1(S)$ is the holonomy of $\rho$,
 \end{itemize}
  such that the operator $\Hh_N^{red}(T_{C_\phi},\tilde b,\bw)^T$   
  is a QHI intertwiner which intertwins $\rho_\lambda$
  and $\rho_{\lambda'}$ as in Proposition \ref{QHIlocrep}. The load of $\rho$ is determined
  by the values of the weight $\kappa$ at the meridian
  curves that form $\bar S \cap \bar M_\phi$.

(2) The family of QHI intertwiners $\{\Hh_N^{red}(T_{C_\phi},\tilde b,\bw)^T\}$
consists of the QT intertwiners which intertwin $\rho_\lambda$ to $\rho_{\lambda'}$ and whose
    traces are well defined invariants of triples $(M_\phi,\rG,\kappa)$
    such that the restriction of $\rG$ to $\pi_1(S)$ is the holonomy of
    $\rho$.
}
\medskip

\Dim (1) Let $\rho\colon\TG^q_{\lambda}\to {\rm End}( (\C^N)^{\otimes 2m})$ be a local representation with parameters $(\{x(e)\}_{e\in \lambda^{(1)}},h)$. Recall Proposition \ref{gluing-var} and the discussion that follows. Assume that the augmented character $[(\rr,\{\xi_\Gamma\}_\Gamma)]$ associated to $\rho$ lies in $i^*(X_0(M_\phi))$, and can be encoded by a system of shape parameters $w\in A_0 \subset A$, the simply connected open neighborhood of $w_h$ in $G(T,\tilde b)$ chosen before Remark \ref{qpar}. Then we can lift $w$ to a system of quantum shape parameters $\bw\in A_{0,N}$, and using the QH-triangulation $(T_{C_\phi},\tilde b,\bw)$, we can define a local representation $\rho'\colon\TG^q_{\lambda}\to {\rm End}( (\C^N)^{\otimes 2m})$ in the same way as $\rho_0$ in the proof of Proposition \ref{QHIlocrep} (see \eqref{relloc2}). It has quantum shear-bend coordinates $\bx(e)$, $e\in \lambda^{(1)}$, and parameters $(\{x(e)\}_{e\in \lambda^{(1)}},h')$, where $\textstyle h'=\prod_{e\in \lambda^{(1)}} \bx(e)$. The system $\{x(e)\}_{e\in \lambda^{(1)}}$ is the same for $\rho$ and $\rho'$, so in order that $\rho'\cong \rho$, it remains to show that we can choose $\bw$ so that $h'=h$. But $\textstyle (h')^N = \prod_{e\in \lambda^{(1)}} x(e) = h^N$, so $h'=\zeta h$ for some $N$-th root of unity $\zeta$. On the other hand, by Theorem \ref{SSum} (1), at all the meridian curves $m_j$ of $\partial \bar M_\phi$ we have
$$\kappa_N(\bw)(m_j)^{N} = \pm a_j^{-2}$$
where $a_j$ is defined before the identity \eqref{geomintpar}, and we note that $a_j^{-2}$ is the dilation factor of the similarity transformation associated to the conjugate of $\rr(m_j)$ fixing $\infty$, whence the squared eigenvalue at $m_j$ selected by $[(\rr,\{\xi_\Gamma\}_\Gamma)]$. Hence 
$$(h')^{2N} = a_1^{-2}\ldots a_r^{-2} = \pm \kappa_N(\bw)(m_1)^{N}\ldots \kappa_N(\bw)(m_r)^{N}.$$
Moreover, there is $\bw'\in A_{0,N}$ so that $\kappa_N(\bw')(m_1)\ldots \kappa_N(\bw')(m_r)$ is any given $N$-th root of $a_1^{-2}\ldots a_r^{-2}$. Let $\zeta'$ be the $N$-th root of unity such that $(h')^{2} = \zeta'\kappa_N(\bw)(m_1)\ldots \kappa_N(\bw)(m_r)$. By using the formula \eqref{cuspweightform} it is easy to check that
\begin{equation}\label{formmj}
\kappa_N(\bw)(m_j) = (-q^{-1})^{t_j}\prod_{e_\vert p_j \in \partial e} \bW^+(e)
\end{equation} 
where $\bW^+(e)$ is the lateral quantum shape parameter of $(S,\lambda)$ at the edge $e$ of $\lambda$, $t_j$ is the number of spikes of the triangles of $\lambda$ at the $j$-th puncture $p_j$ of $S$, and the product is over all the edges $e$, counted with multiplicity, having $p_j$ as an endpoint. Recalling that there are $3m=-3\chi(S)$ edges in $\lambda$, we deduce
\begin{align}
\kappa_N(\bw)(m_1)\ldots \kappa_N(\bw)(m_r) & =  (-q^{-1})^{3.2m} \prod_{e\in \lambda^{(1)}} \bW^+(e)^2 \label{lastform0}\\
 & =  (-q^{-1})^{6m}(-q)^{6m} \prod_{e\in \lambda^{(1)}} \bx(e)^2 =  (h')^2.\label{lastform1}
 \end{align}
Hence $\zeta'=1$. Then take $\bw'\in A_{0,N}$ so that $\kappa_N(\bw')(m_1)\ldots \kappa_N(\bw')(m_r) = (h')^{2}\zeta^{-2}$. The load $h''$ of the representation $\rho''$ associated to $\bw'$ satisfies $(h'')^2 = (h')^{2}\zeta^{-2} = h^2$. Since $N$ is odd, and again $(h'')^N=h^N$, eventually $h=h''$.
This achieves the proof of the claim (1) of the theorem.
\smallskip

(2) We need to describe the general form of the QT intertwiners associated to sequences of diagonal exchanges $\lambda\rightarrow \ldots \rightarrow \lambda'$, for arbitrary standard local representations $\rho_{\lambda}$, $\rho_{\lambda'}$ of $\TG_{\lambda}^q$, $ \TG_{\lambda'}^q$. This is done in Lemma \ref{allQT}. For that, we reconsider the arguments of the proof of the claim (1) in the First Main Theorem. First, we look at the case of a diagonal exchange $\lambda \rightarrow \lambda'$, occuring in squares $Q$, $Q'$ as in Figure \ref{typefig3}. We fix standard local representations $\rho_{\lambda}$, $\rho_{\lambda'}$ of $\TG_{\lambda}^q$, $ \TG_{\lambda'}^q$ such that $\rho_{\lambda'}\circ (\varphi_{\lambda\lambda'}^q)^{-1}$ is isomorphic to $\rho_{\lambda}$, and differs from it only on $Q$, $Q'$, where $\rho_{\lambda}$, $\rho_{\lambda'}$ are represented by irreducible representations $\rho_j$ of the triangle algebra with parameters $\by^j_k$, where $j\in \{0,1,2,3\}$ and $k\in\{0,1,2\}$. 
\begin{lem}  \label{roughw} With the above notations we have
\begin{equation}\label{relexdiag0}
\by^3_0/\by^2_0 = \by^1_2/\by^0_2\ , \ \by^1_1/\by^2_1 = \by^3_2/\by^0_0\ ,\ \by_1^3\by_0^1 = (\by^2_2\by^0_1)^{-1}.
\end{equation}
Moreover, denoting these scalars by $\tilde \bw_0$, $\tilde \bw_1$ and $\tilde \bw_2$ respectively, $(w_0,w_1,w_2) := (\tilde\bw_0^N,\tilde\bw_1^N,\tilde\bw_2^N)$ is a triple of shape parameters on $(\Delta,b)$, that is $w_{i+1} = (1-w_i)^{-1}$ for $i=0,1$.
\end{lem}
We say that $(\tilde \bw_0,\tilde\bw_1,\tilde\bw_2)$ is a triple of \emph{rough} $q$-shape parameters on $(\Delta,b)$. Each of the scalars $\tilde \bw_j$ corresponds to a pair of opposite edges of $\Delta$. Note that, in general, the product $\tilde \bw_0\tilde \bw_1\tilde \bw_2$ can be any $N$-th root of $-1$, the special case of $-q$ being achieved by the $q$-shape parameters.  
\medskip

\Dim We prove first the second claim, that is, the $N$-th powers of the scalars set equal in \eqref{relexdiag0} are indeed equal, and form a triple of shape parameters $(w_0,w_1,w_2)$. By the relations \eqref{relq} we have
$$\begin{array}{ll}
\varphi_{\lambda\lambda'}^q((X_5')^N) = \varphi_{\lambda\lambda'}^q((X_5'))^N = X_5^{-N}\\
\varphi_{\lambda\lambda'}^q((X_1')^N) = \varphi_{\lambda\lambda'}^q(X_1')^N = (X_1+qX_5X_1)^N = X_1^N+X_5^NX_1^N\\
\varphi_{\lambda\lambda'}^q((X_2')^N) = \varphi_{\lambda\lambda'}^q(X_2')^N = (X_2^{-1}+qX_2^{-1}X_5^{-1})^{-N} = (X_2^{-N}+X_2^{-N}X_5^{-N})^{-1}
\end{array}$$
and a relation like the second for $X_3'$, $X_3$, and like the fourth for $X_4'$, $X_4$. Note that we use the $q$-binomial formula and the fact that $q^2$ is a primitive $N$-th root of $1$, $N$ odd, to deduce the results. Hence the shear-bend parameters $x_i$, $x_i'$ at the edges $e_i$, $i\in \{1,\ldots,5\}$, satisfy
$$x_5' = x_5^{-1}\ ,\ x_1' = x_1(1+x_5)\ ,\ x_3' = x_3(1+x_5)\ ,\ x_2' = x_2(1+x_5^{-1})^{-1}\ ,\ x_4' = x_4(1+x_5^{-1})^{-1}.$$
Set $x_5 = -w_2$. By a mere rewriting of formulas, it follows easily from this, \eqref{relloc} and $\bx(e)^N=x(e)$ for all edges $e$ of $\lambda$, that 
\begin{equation}\label{2relexdiag}
\begin{array}{c}
- w_2 = (\by_1^3\by_0^1)^N = (\by^2_2\by^0_1)^{-N}\\
w_0=x_1/x_1' = x_3/x_3' = (\by^3_0/\by^2_0)^N = (\by^1_2/\by^0_2)^N\\
w_1 = x_2/x_2' = x_4/x_4' = (\by^1_1/\by^2_1)^N = (\by^3_2/\by^0_0)^N.
\end{array}
\end{equation}
and $w_{i+1} = (1-w_i)^{-1}$ for $i=0,1$. Now we prove \eqref{relexdiag0}. Let us reconsider Proposition \ref{Sintertwins} in the present situation. Instead of $\rho_\lambda$ and the system $(\bw_0,\bw_1,\bw_2)$ of quantum shape parameters on $(\Delta,b)$, it is the pair of local representations $\rho_\lambda$, $\rho_{\lambda'}$ that are given on $Q$, $Q'$, and we are looking for $(\tilde\bw_0,\tilde\bw_1,\tilde\bw_2)$ such that the identity \eqref{intprop} holds true. Again \eqref{intprop} yields the first equality of \eqref{firstrel}, which is the last relation of \eqref{relexdiag0}, without implying any constraints on $(\tilde\bw_0,\tilde\bw_1,\tilde\bw_2)$. Then, use it to define $\tilde\bw_2$ by the same formula. Replacing all uses of the relation \eqref{cycdil} by \eqref{cycdil2} in the other computations defines $\tilde\bw_0 = \by^3_0/\by^2_0$ and $\tilde\bw_1 = \by^1_1/\by^2_1$, and eventually proves the first two identities in \eqref{relexdiag0}. For instance, in the case of $X=X_4$ the result becomes
$$\Ll_N^T(\Delta,b,\tilde \bw) \circ (\by^3_2 A_2 \otimes I_N)= (1-q^2\tilde\bw_2^{-1} A_2 \otimes A_1) \circ \left(I_N \otimes (-q^{-1}\tilde\bw_0\tilde\bw_2\by^3_2 A_0)\right) \circ \Ll_N^T(\Delta,b,\tilde\bw)$$
and hence \eqref{intprop4} holds true whenever $\tilde\bw_0 = -q\by^0_0/\by^3_2\tilde\bw_2$. Comparing with the result of the similar computation for $X=X_2$ we find  the second identity in \eqref{relexdiag0}. The other cases are similar. \hfill$\Box$
\medskip

Now we have:
\begin{prop} {\rm (See \cite{B-L,Filippo})} Assume that $S$ is an ideal polygon with $p\geq 3$ vertices. Then for every ideal triangulation $\lambda$ of $S$, every local representation of $\TG_\lambda^q$ is irreducible.
\end{prop}
In the situation of the lemma we have the ideal squares $Q$, $Q'$, and so the proposition implies that $\Ll_N^T(\Delta,b,\tilde \bw)$ is a representative of the unique projective class of intertwiners from $\rho_\lambda$ to $\rho_{\lambda'}\circ (\varphi_{\lambda\lambda'}^q)^{-1}$. More generally, let $\lambda$, $\lambda'$ be ideal triangulations of $S$ connected by a sequence of diagonal exchanges, and let $\rho$, $\rho'$ be isomorphic standard local representations of $\TG_S^q$. Then every QT intertwiner $L^{\rho\rho'}_{\lambda\lambda'}$ between $\rho_\lambda$ and $\rho_{\lambda'}'$ has a representative of the form \eqref{decompositionint}, where $L^{\rho\rho'}_{\lambda_i\lambda_{i+1}} =\Ll_N^T(\Delta,b,\tilde \bw)$ for some system $\tilde \bw:=(\tilde \bw_0,\tilde\bw_1,\tilde\bw_2)$ of rough 
$q$-shape parameters (as in Lemma \ref{roughw}) on the tetrahedron $(\Delta,b)$ associated to the $i$-th diagonal exchange. 

Hence, proceeding as in the proof of Proposition \ref{QHIlocrep} we deduce that:
\begin{lem}\label{allQT} Every projective class of QT intertwiner in $\mathcal{L}_{\lambda\lambda'}^{\rho\rho'}$ is represented by a QH state sum $\Hh_N^{red}(T_{C_\phi},\tilde b,\tilde \bw)^T$, where $\tilde\bw$ is a system of  rough $q$-shape parameters on $T_{C_\phi}$.  
\end{lem}
The proof of the claim (2) of the Second Main Theorem then follows by analysing the invariance properties of the QH state sums ${\rm Trace}\left(\Hh_N^{red}(T_{C_\phi},\tilde b,\tilde \bw)\right)$. Such an analysis has been done in \cite{GT} and \cite{AGT} in the case of the unreduced QH state sums on arbitrary closed pseudo-manifold triangulations. It applies verbatim to the reduced QH state sums by combining the results of Section 5 of \cite{GT} with Proposition 6.3 (1) and 8.3 of \cite{NA}. 

The result is that ${\rm Trace}\left(\Hh_N^{red}(T_{C_\phi},\tilde b,\tilde \bw)\right)$ is an invariant of the triple $(M_\phi,\rG,\kappa)$ if and only if it is invariant under any layered `$2\leftrightarrow 3$ transit'  or `lune transit' of $(T_{C_\phi},\tilde b,\tilde \bw)$; such transits enhance to $\tilde b$ and $\tilde \bw$ the usual $2\leftrightarrow 3$ and lune moves between layered triangulations (which correspond to the pentagon and square relations for diagonal exchanges between surface triangulations). This invariance property happens if and only if $\tilde \bw$ satisfies the tetrahedral and edge relations along every interior edge of $T_{C_\phi}$, that is, when $\tilde \bw$ is a genuine system of quantum shape parameters. In that situation ${\rm Trace}\left(\Hh_N^{red}(T_{C_\phi},\tilde b,\tilde \bw)\right)$ a QHI intertwiner.

For instance, the relation \eqref{firstrel} implies that, modifying $T_{C_\phi}$ by introducing in the sequence $\lambda \rightarrow \ldots \rightarrow \lambda'$ two consecutive diagonal exchanges along an edge $e$, $\Hh_N^{red}(T_{C_\phi},\tilde b,\tilde \bw)$ is unchanged if and only if the total quantum shape parameter at $e$ satisfies $\bW(e)=q^2$.  
\medskip

This concludes the proof of our Second Main Theorem.\cvd
\medskip

Finally, let us consider Corollary \ref{cor0} and \ref{cor}. For the reader's convenience we re-state them below as Corollary 1 and 2. We keep the usual setting. Let us fix a QH-triangulation $(T_{C_\phi},\tilde b,\bw)$ of the cylinder $C_\phi$. Recall that $\Hh_N^{red}(M_{\phi},\rG,\kappa) = {\rm Trace}\big(\Hh_N^{red}\left(T_{C_\phi},\tilde b,\bw)\right)$, and that $M_\phi$ is the interior of a compact manifold $\bar M_\phi$ with boundary made by tori. Call {\it longitude} any simple closed curve in $\partial \bar M_\phi$ intersecting a fibre of $\bar M_\phi$ in exactly one point. 
\medskip

\noindent {\bf Corollary 1.} The reduced QH invariants $\Hh_N^{red}(M_{\phi},\rG,\kappa)$ do not depend on the values of the weight $\kappa$ on the longitudes.
\medskip

\Dim Denote by $\rho_\lambda$ the local representation of $\TG_\lambda^q$ associated to $(T_{C_\phi},\tilde b,\bw)$ as in Proposition \ref{QHIlocrep}, where $\lambda=\lambda_0$. Since $\kappa = \kappa_N(\bw)$ is given by its values on the longitudes and the meridian curves $m_j$ of $\partial \bar{M}_\phi$, it is enough to show that $\rG$ and the scalars $\kappa(m_j)$ determine the parameters $(\{x_1,\ldots,x_n\},h)$ of $\rho_\lambda$. We already know that $x_1,\ldots,x_n$ are determined by $\rG$. Now, by the identities \eqref{lastform0} and \eqref{lastform1} we have $h^2 = \kappa(m_1)\ldots \kappa(m_r)$. Among the two square roots of $h^2$, the load $h$ is the only one which is an $N$-th root of $(-1)^r a_1^{-1}\ldots a_r^{-1}$ (using the notations of \eqref{geomintpar}). The eigenvalues $a_i$ being determined by the augmented character $\rG$, the conclusion follows. \cvd
\medskip

Recall that given a local representation $\rho_\lambda$ and an irreducible representation $\mu$ of $\TG_{\lambda}^q$, we denote by $\rho_\lambda(\mu)$ the isotypical component associated to $\mu$ in the direct sum decomposition of $\rho_\lambda$ into irreducible factors. Consider the {\it isotypical intertwiners} $L^\phi_{\rho_\lambda(\mu)} := \Hh_N^{red}(T_{C_\phi},\tilde b',\bw)^T_{\vert \rho_\lambda(\mu)}$.
\medskip

\noindent {\bf Corollary 2.} (1) The trace of $L^\phi_{\rho_\lambda(\mu)}$ is an invariant of $(M_{\phi},\rG,\kappa)$ and $\mu$, well-defined up to multiplication by $4N$-th roots of unity. It depends on the isotopy class of $\phi$ and satisfies
  $$\Hh_N^{red}(M_{\phi},\rG,\kappa)  = \sum_{\rho_\lambda(\mu)\subset \rho_\lambda}{\rm Trace}\left(L^\phi_{\rho_\lambda(\mu)}\right).$$
\noindent (2) The invariants ${\rm Trace}\left(L^\phi_{\rho_\lambda(\mu)}\right)$ do not depend on the values of the weight $\kappa$ on longitudes. 
\medskip

\Dim (1) We take back the notations used for Corollary 1. Once again, $\Hh_N^{red}(M_{\phi},\rG,\kappa)$ is determined by $\phi$ and the isomorphism class of the local representation $\rho_\lambda$, because it is the trace of $\Hh_N^{red}(T_{C_\phi},\tilde b,\bw)$. Each isotypical component of $\rho_\lambda$ is the intersection of one eigenspace of each of the so called {\it puncture elements} of $\TG_\lambda^q$ (see \cite{B-L} and \cite{Tou}). Hence it is determined by the parameters $(\{x_1,\ldots,x_n\},h)$ of $\rho_\lambda$, together with a system of {\it puncture weights} $p_j\in \mc^*$, the eigenvalues of the puncture elements. These can be any $N$-th root of $\kappa(m_j)^N$, where $\kappa(m_j)$ is computed in \eqref{formmj}. Hence $\rG$ and $\kappa$ determine one isotypical component of $\rho_\lambda$. The others correspond to all other systems of puncture weights obtained by multiplying each $p_j$ with some $N$-th root of unity. Since $\Hh_N^{red}(T_{C_\phi},\tilde b,\bw)$ intertwins $\rho_\lambda$ and $\rho_{\lambda'}$, it intertwins their isotypical components too. The usual invariance proof of the reduced QHI still apply, and we get that the traces of the isotypical intertwiners are singly invariant.

By the results of \cite{AGT,NA}, the invariants $\Hh_N^{red}(M_{\phi},\rG,\kappa)$ actually depend on the choice of mapping torus realization $M_\phi$ of $M$. Therefore the invariants ${\rm Trace}(L^\phi_{\rho_\lambda(\mu)})$ do as well. By the same arguments as for $\Hh_N^{red}(M_{\phi},\rG,\kappa)$, they are also defined for any augmented character $\rG$ in a Zariski open subset of the eigenvalue subvariety of $X_0(M)$ (see the comments on Theorem \ref{MAINTEO2} in the Introduction). 
 
The proof of (2) is the same as the one of the previous corollary. \cvd

\end{document}